\documentclass[11pt]{amsart}

\usepackage{lahn}
\usepackage{ALLM}
\usepackage{fullpage}

\title{A note on an effective characterization of covers with an application to higher rank representations}

\author{Tarik Aougab}
\address{Department of Mathematics \\ Haverford College \\ Haverford, PA 19041}
\email{taougab@haverford.edu}

\author{Max Lahn}
\address{Department of Mathematics \\ University of Michigan \\ Ann Arbor, MI 48109}
\email{maxlahn@umich.edu}

\author{Marissa Loving}
\address{Department of Mathematics \\ University of Wisconsin \\ Madison, WI 53706}
\email{mloving2@wisc.edu}

\author{Nicholas Miller}
\address{Department of Mathematics \\ University of Oklahoma \\ Norman, OK 73019}
\email{nickmbmiller@ou.edu}

\begin{document}

\begin{abstract}


In this note we prove an effective characterization of when two finite-degree covers of a connected, orientable surface of negative Euler characteristic are isomorphic in terms of which curves have simple elevations, weakening the hypotheses to consider curves with explicitly bounded self-intersection number. As an application we show that for sufficiently large $ N $, the set of unmarked traces associated to simple closed curves in a generically chosen representation to $ \SL_{ N } \of{ \R } $ distinguishes between pairs of non-isomorphic covers. 


\end{abstract}

\maketitle

\section{Introduction}

In the setting of higher Teichm\"{u}ller theory, marked spectral rigidity results are known for some of the many class functions on $ \pi_{ 1 } \of{ S } $ which share roles analogous to hyperbolic length. For example, Bridgeman--Canary--Labourie \cite{BCL20} show that a Hitchin representation of $ \pi_{ 1 } \of{ S } $ is determined by either the marked spectral radii of simple curves or the marked traces associated to simple curves. In this setting, we show that covers can be distinguished by the set of unmarked traces of simple closed curves.

\begin{restatable}{theorem+}{SLnTrace} \label{thm:SLnTrace}

There is a positive integer $ N $, explicitly computable in terms of $ \abs{ \chi \of{ S } } $, $ \deg \of{ p } $, and $ \deg \of{ q } $, so that $ p $ and $ q $ are isomorphic if and only if the restrictions of linear representations $ \pi_{ 1 } \of{ S } \to \SL_{ N } \of{ \R } $ or $ \pi_{ 1 } \of{ S } \to \SL_{ N } \of{ \C } $ to $ \pi_{ 1 } \of{ X } $ and $ \pi_{ 1 } \of{ Y } $ are generically simple trace isospectral.

\end{restatable}

We call the reader's attention to the words \emph{explicitly computable} in the above statement of \autoref{thm:SLnTrace}. We obtain this by first proving the following effective version of \cite[Theorem~3]{ALLM23}, with its hypotheses weakened so that we need only consider the elevations of curves on $ S $ with explicitly bounded self-intersection number.

\begin{restatable}{theorem+}{effective} \label{thm:effective}

There is an integer $ M $, explicitly computable in terms of $ \abs{ \chi \of{ S } } $, $ \deg \of{ p } $, and $ \deg \of{ q } $, so that $ p $ and $ q $ are isomorphic as covers of $ S $ if and only if for every curve $ \gamma $ on $ S $ with $ i \of{ \gamma , \gamma } \leq M $, $ \gamma $ has a simple elevation along $ p $ to $ X $ if and only if it has a simple elevation along $ q $ to $ Y $.

\end{restatable}

We note that although a non-effective version of \autoref{thm:effective} can be derived quickly from \cite[Theorem~3]{ALLM23} (such an argument is outlined at the beginning of \autoref{sec:effective}), justifying that the integer $ M $ is \emph{explicitly computable} in terms of $ \abs{ \chi \of{ S } } $, $ \deg \of{ p } $, and $ \deg \of{ q } $ requires an entirely different argument than that used in the proof of \cite[Theorem~3]{ALLM23}. In particular, \autoref{thm:effective} gives a completely distinct proof of \cite[Theorem~3]{ALLM23}. We decided to relegate this result to its own note rather than including it in \cite{ALLM23} because achieving effective control in \autoref{thm:effective} requires navigating through some computationally heavy lemmas; we use fine properties of the geometry of the curve complex to facilitate these arguments. Our proof also requires effective versions of the main theorem of \cite{Riv01} and of \cite[Proposition~6.2]{Tan17}, which we prove in detail in \autoref{sec:Rivin} and \autoref{sec:effectivizingTang}, respectively, as they may be of independent interest.

\subsection{Organization of paper}
\label{sub:organization}

We provide some background definitions and tools in \autoref{sec:preliminaries}. In \autoref{sec:effective}, we prove \autoref{thm:effective}. We then prove \autoref{thm:SLnTrace} in \autoref{sec:SLnTrace}. Finally, \autoref{sec:Rivin} and \autoref{sec:effectivizingTang} contain the proofs of effective versions of the main theorem of \cite{Riv01} and \cite[Proposition~6.2]{Tan17}, respectively, which are needed for the proof of \autoref{thm:effective}.

\subsection{Acknowledgements}

We gratefully acknowledge support from NSF grants DMS-1939936 (Aougab), DMS-1906441 (Lahn), DMS-1902729 \& DMS-2231286 (Loving), and DMS-2005438/2300370 (Miller).


\section{Preliminaries}\label{sec:preliminaries} 

The interested reader can refer to \cite{ALLM23} for a more in depth treatment of much of the relevant background. To keep this note concise, we will only mention two of the most specialized results that we make use of here which are both related to the geometry of the curve complex.

\subsection{Geometry of the curve complex}\label{sec:curvegraph}

We will make use of the following work of Aougab-Taylor, which constructs geodesic rays whose vertices have controlled intersections.

\begin{theorem*}[{\cite[Theorem~1.1]{AT14}}] \label{thm:aougab-taylor}

For any simple curve $ \eta \in \CC \of{ S } $ on a closed surface $ S $ of negative Euler characteristic, there is a sequence $ \of{ v_{ k } }_{ k = 0 }^{ \infty } $ in $ \CC \of{ S } $ with $ v_{ 0 } = \eta $, $ d_{ S } \of{ v_{ j } , v_{ k } } = \abs{ j - k } $ for all indices $ j , k \in \N $, and
\[
i \of{ v_{ j } , v_{ k } } \leq \of{ D + 3 }^{ 2 \max \of{ j , k } - 5 } \of{ \frac{ 3 \abs{ \chi \of{ S } } }{ 2 } }^{ \abs{ j - k } + 2 } + f_{ j , k } \of{ \abs{ \chi \of{ S } } } ,
\]
where $ D $ is a fixed, universal constant described in \autoref{thm:boundedGeodesicImage}, and $ f_{ j , k } \of{ \abs{ \chi \of{ S } } } $ is $ O \of{ \abs{ \chi \of{ S } }^{ \abs{ j - k } - 4 } } $.

\end{theorem*}

\noindent We remark that the above bound on intersection number takes a slightly different form than the one found in \cite{AT14} since we are only considering closed surfaces.

As proven by Rafi--Schleimer \cite{RS09}, finite-degree covers of surfaces give rise to quasi-isometric embeddings of their curve complexes. Specifically, given a finite-degree cover $ p \colon X \to S $, we define $ \widetilde{ p } \colon \CC \of{ S } \to \CC \of{ X } $ so that for each simple curve $ \gamma \in \CC \of{ S } $, $ \widetilde{ p } \of{ \gamma } $ is a choice of elevation of $ \gamma $ along $ p $ to $ X $. At first glance, our definition of $ \widetilde{ p } $ appears to rely significantly on a choice of elevation for every simple curve on $ S $ and thus its properties will not be well behaved when defined as described. However, any two elevations of a simple curve $ \gamma $ along $ p $ are disjoint on $ X $, since any intersection between distinct elevations on $ X $ would descend to a self-intersection of $ \gamma $ on $ S $. Thus different choices of elevations for $ \widetilde{ p } $ yield maps which have uniformly bounded distance at most $ 1 $, and $ \widetilde{ p } \colon \CC \of{ S } \to \CC \of{ X } $ is said to be \emph{coarsely well-defined}.

\begin{theorem*}[{\cite[Theorem~7.1]{RS09}}] \label{thm:CCembedding}

If $ p \colon X \to S $ is a finite-degree cover of a closed, orientable surface $ S $ of negative Euler characteristic, then $ \widetilde{ p } \colon \CC \of{ S } \to \CC \of{ X } $ is a quasi-isometric embedding.

\end{theorem*}

Note that in \cite{RS09}, the above theorem is formulated in terms of the many-to-one relation which associates a simple curve on the base surface to its full pre-image on the cover, which is a multi-curve. Since a multi-curve has diameter at most $ 1 $ in the curve complex, this distinction will not matter for our purposes.

We also note that the quality of the above quasi-isometry induced by such a finite-degree cover of surfaces can be explicitly bounded in terms of the Euler characteristic of the base and the degree of the cover. We require such bounds, as computed by Aougab--Patel--Taylor \cite[Theorem~7.1]{APT22}, in \autoref{sec:effective}.

Since the image $ \widetilde{ p } \of{ \CC \of{ S } } $ of $ \widetilde{ p } \colon \CC \of{ S } \into \CC \of{ X } $ is quasi-convex in the Gromov hyperbolic space $ \CC \of{ X } $, nearest point projection $ \nu_{ S } \colon \CC \of{ X } \onto \widetilde{ p } \of{ \CC \of{ S } } $ is coarsely well-defined and distance non-increasing. We will make use of the following result of Tang.

\begin{restatable}[{\cite[Proposition~6.2]{Tan17}}]{theorem*}{Tang} \label{thm:Tang}

Let $ p \colon X \to S $ be a finite-degree regular cover of a closed surface $ S $ of negative Euler characteristic. There is a uniform upper bound, depending only on $ \abs{ \chi \of{ S } } $ and $ \deg \of{ p } $, on the distance between any circumcenter of the $ \Deck \of{ p } $-orbit of a simple curve $ \alpha \in \CC \of{ X } $ and its nearest point projection $ \nu_{ S } \of{ \alpha } $.

\end{restatable}

\noindent See \autoref{sec:effectivizingTang} for an explicit computation of this uniform upper bound.

\section{Elevations of curves with bounded self-intersection number} \label{sec:effective}

In this section, we prove \autoref{thm:effective}, which we restate below for the reader's convenience.

\effective*

We stress the importance of the words \textit{explicitly computable} in the above theorem statement. Indeed, the fact that some positive integer $ M $ satisfies the conclusions of \autoref{thm:effective} follows immediately from \cite[Theorem~3]{ALLM23} and the following argument:

Let $ d = \max \of{ \deg \of{ p } , \deg \of{ q } } $, and note that there are only finitely many isomorphism classes of covers of $ S $ of degree at most $ d $. In particular, there are only finitely many \emph{pairs} of distinct isomorphism classes of covers of $ S $ of degree at most $ d $. Choose a representative of each isomorphism class and note that, for each pair of non-isomorphic covers in this finite list, there is some curve $ \gamma $ on $ S $ which has a simple elevation along one cover but not along the other. Any positive integer $ M $ at least as large as the largest self-intersection number of these curves $ \gamma $ will satisfy the conclusions of \autoref{thm:effective}.

Thus, in the following proof sketch, we take care to justify that $ M $ is, at least in principle, explicitly computable from the Euler characteristic $ \chi \of{ S } $ and the degrees $ \deg \of{ p } $ and $ \deg \of{ q } $. We will use the phrase \emph{explicitly computable} to mean precisely this, and we will call a number \emph{explicitly bounded} if it can be bounded in terms of an explicitly computable number.

For example, the quality of the quasi-isometric embedding of curve complexes induced by a finite-degree cover (for a detailed discussion of this quasi-isometric embedding, see \autoref{thm:CCembedding}) is explicitly bounded in terms of the degree of the cover and the Euler characteristic of the base surface. This is the content of the following result of Aougab-Patel-Taylor.

\begin{theorem*}[{\cite[Theorem~7.1]{APT22}}] \label{thm:APT22}

If $ p \colon X \to S $ is a finite-degree cover of a closed, orientable surface $ S $ of negative Euler characteristic, then
\[
\frac{ d_{ S } \of{ \alpha , \beta } }{ 80 e^{ 54 } \pi \deg \of{ p } \cdot \abs{ \chi \of{ S } }^{ 13 } } \leq d_{ X } \of{ \widetilde{ p } \of{ \alpha } , \widetilde{ p } \of{ \beta } } \leq d_{ S } \of{ \alpha , \beta },
\]
for all simple curves $ \alpha , \beta \in \CC \of{ S } $.

\end{theorem*}

\noindent We will require this effective statement for the proof of \autoref{thm:effective}.

\subsection{Counting simple curves of bounded hyperbolic length}

The following theorem of Rivin gives effective estimates on the growth of the number of simple geodesics of bounded length on a hyperbolic surface of finite type as a polynomial in the bound. The degree of this polynomial increases with the complexity of the surface, so that the number of simple geodesics of bounded length grows noticeably faster on surfaces of higher genus. By a combinatorial argument, we show that covers of different degrees can be distinguished by which curves of explicitly bounded length have simple elevations.

\begin{restatable}[\cite{Riv01}]{theorem*}{Rivin} \label{thm:Rivin}

For any hyperbolic surface $ \Sigma $ of finite type $ \of{ g , b , n } $, there are positive real numbers $ c_{ 1 } \of{ \Sigma } $, $ c_{ 2 } \of{ \Sigma } $, and $ L_{ 0 } \of{ \Sigma } $ so that
\[
c_{ 1 } \of{ \Sigma } \leq \frac{ \mathcal{ N } \of{ L , \Sigma } }{ L^{ 6 g - 6 + 2 b + 2 c } } \leq c_{ 2 } \of{ \Sigma },
\]
for all real numbers $ L \geq L_{ 0 } \of{ \Sigma } $.

\end{restatable}

Specifically, we will use the following technical corollary, which utilizes the reliance of the above bounds on the genus $ g $ to find curves of bounded length which have simple elevations to some covers but not to others.

\begin{corollary*} \label{cor:awfulCombinatorics}

Consider a finite collection $ \set{ \pi_{ k } }_{ k = 0 }^{ n } $ of finite-degree covers $ \pi_{ k } \colon \Sigma_{ k } \to \Sigma $ of a closed, hyperbolic surface $ \Sigma $. If $ \deg \of{ \pi_{ 0 } } > \deg \of{ \pi_{ k } } $ for all indices $ 1 \leq k \leq n $, then there is a simple closed geodesic $ \alpha $ on $ \Sigma_{ 0 } $ whose projection $ \pi_{ 0 } \of{ \alpha } $ has no simple elevations along $ \pi_{ k } $ to $ \Sigma_{ k } $ for all indices $ 1 \leq k \leq n $.

Moreover, such a simple closed geodesic $ \alpha $ can be chosen so that $ \len_{ \Sigma_{ 0 } } \of{ \alpha } $ is explicitly bounded in terms of $ n $, the degree $ \deg \of{ \pi_{ 0 } } $, the Euler characteristic $ \abs{ \chi \of{ \Sigma } } $, and the constants guaranteed by \autoref{thm:Rivin}.

\begin{proof}[Proof of \autoref{cor:awfulCombinatorics}]

For each index $ 0 \leq k \leq n $, let $ c_{ 1 } \of{ \Sigma_{ k } } $, $ c_{ 2 } \of{ \Sigma_{ k } } $, and $ L_{ 0 } \of{ \Sigma_{ k } } $ denote the constants guaranteed by \autoref{thm:Rivin}. 

\begin{claim} \label{clm:RivinComputation}

There is a positive real number $ L $, explicitly computable in terms of $ n $, the degree $ \deg \of{ \pi_{ 0 } } $, the Euler characteristic $ \abs{ \chi \of{ \Sigma } } $, and the constants guaranteed by \autoref{thm:Rivin}, so that
\[
\frac{ \mathcal{ N } \of{ L , \Sigma_{ 0 } } }{ \deg \of{ \pi_{ 0 } } } > n \max_{ 1 \leq k \leq n }{ \mathcal{ N } \of{ \deg \of{ \pi_{ k } } L , \Sigma_{ k } } } .
\]

\end{claim}

Although we include the following computation for completeness, the key point is that $L$ is explicitly computable -- its exact value plays little role in the overall argument. As such, we encourage skipping the following proof of \autoref{clm:RivinComputation} on a first reading.

\begin{proof}[Proof of \autoref{clm:RivinComputation}]

If $ L \geq L_{ 0 } \of{ \Sigma_{ k } } $ for all indices $ 0 \leq k \leq n $, \autoref{thm:Rivin} implies that
\begin{align*}
n \mathcal{ N } \of{ \deg \of{ \pi_{ k } } L , \Sigma_{ k } } & \leq n c_{ 2 } \of{ \Sigma_{ k } } \of{ \deg \of{ \pi_{ k } } L }^{ 3 \abs{ \chi \of{ \Sigma_{ k } } } } \\
& \leq n c_{ 2 } \of{ \Sigma_{ k } } \deg \of{ \pi_{ 0 } }^{ 3 \deg \of{ \pi_{ 0 } } \abs{ \chi \of{ \Sigma } } } L^{ 3 \deg \of{ \pi_{ k } } \abs{ \chi \of{ \Sigma } } } ,
\end{align*}
for all indices $ 1 \leq k \leq n $, and similarly that
\begin{align*}
\frac{ \mathcal{ N } \of{ L , \Sigma_{ 0 } } }{ \deg \of{ \pi_{ 0 } } } \geq \frac{ c_{ 1 } \of{ \Sigma_{ 0 } } L^{ 3 \abs{ \chi \of{ \Sigma_{ 0 } } } } }{ \deg \of{ \pi_{ 0 } } } = \frac{ c_{ 1 } \of{ \Sigma_{ 0 } } L^{ 3 \deg \of{ \pi_{ 0 } } \abs{ \chi \of{ \Sigma } } } }{ \deg \of{ \pi_{ 0 } } } .
\end{align*}
It therefore suffices to find an explicitly computable solution $ L $ for the inequalities
\[
\frac{ c_{ 1 } \of{ \Sigma_{ 0 } } L^{ 3 \deg \of{ \pi_{ 0 } } \abs{ \chi \of{ \Sigma } } } }{ \deg \of{ \pi_{ 0 } } } > n c_{ 2 } \of{ \Sigma_{ k } } \deg \of{ \pi_{ 0 } }^{ 3 \deg \of{ \pi_{ 0 } } \abs{ \chi \of{ \Sigma } } } L^{ 3 \deg \of{ \pi_{ k } } \abs{ \chi \of{ \Sigma } } } ,
\]
for all $k\in\{1,\dots,n\}$.
One can check that
\[
L \coloneqq \max \of{ \max_{ 0 \leq k \leq n }{ L_{ 0 } \of{ \Sigma_{ k } } } , 1 + \max_{ 1 \leq k \leq n }{ \frac{ n c_{ 2 } \of{ \Sigma_{ k } } \deg \of{ \pi_{ 0 } }^{ 1 + 3 \deg \of{ \pi_{ 0 } } \abs{ \chi \of{ \Sigma } } } }{ c_{ 1 } \of{ \Sigma_{ 0 } } } } },
\]
is one such solution. \qedhere

\end{proof}

Let $ U_{ 0 } $ denote the set of simple closed geodesics on $ \Sigma_{ 0 } $ of length at most $ L $, so that $ \card{ U_{ 0 } } = \mathcal{ N } \of{ L , \Sigma_{ 0 } } $. Moreover, denote by $ U \coloneqq \pi_{ 0 } \of{ U_0 } $ the set of closed geodesics on $ \Sigma $ which have elevations along $ \pi_{ 0 } $ to $ \Sigma_{ 0 } $ lying in $ U_{ 0 } $, and note that every curve in $ U $ has length at most $ L $, although such a curve need not be simple. Finally, for each index $ 1 \leq k \leq n $, denote by $ U_{ k } \coloneqq \pi_{ k }^{ -1 } \of{ U } $ the set of (not necessarily simple) elevations along $ \pi_{ k } $ to $ \Sigma_{ k } $ of curves in $ U $, each of which has length at most $ \deg \of{ \pi_{ k } } L $.

For each index $ 0 \leq k \leq n $, the cover $ \pi_{ k } \colon \Sigma_{ k } \to \Sigma $ gives a surjective map $ U_{ k } \onto U $ which is at most $ \deg \of{ \pi_{ k } } $-to-one. Thus
\begin{align*}
\card{ U } & \geq \frac{ \card{ U_{ 0 } } }{ \deg \of{ \pi_{ 0 } } } = \frac{ \mathcal{ N } \of{ L , \Sigma_{ 0 } } }{ \deg \of{ \pi_{ 0 } } } > n \max_{ 1 \leq k \leq n }{ \mathcal{ N } \of{ \deg \of{ \pi_{ k } } L , \Sigma_{ k } } } \geq \sum_{ k = 1 }^{ n }{ \mathcal{ N } \of{ \deg \of{ \pi_{ k } } L , \Sigma_{ k } } } ,
\end{align*}
which implies that there is a closed geodesic in $ U $ for which every corresponding curve in $ \bigcup_{ k = 1 }^{ n }{ U_{ k } } $ is non-simple; that is, there is a simple closed geodesic $ \alpha $ on $ \Sigma_{ 0 } $ of length $ \len_{ \Sigma_{ 0 } } \of{ \alpha } \leq L $ so that for all indices $ 1 \leq k \leq n $, $ \pi_{ 0 } \of{ \alpha } $ has no simple elevations along $ \pi_{ k } $ to $ \Sigma_{ k } $. \qedhere

\end{proof} \setcounter{claim}{0}

\end{corollary*}

See \autoref{sec:Rivin} for a complete and self-contained proof of an effective version of \autoref{thm:Rivin}. Specifically, we will exhibit a dependence of the constants in question on the injectivity radius to choose a metric in the proof of \autoref{thm:effective} so that the constants are explicitly bounded.

\subsection{The collar lemma}

The collar lemma is a classical result relating the geometry of certain regular neighborhoods of a simple closed geodesic on a hyperbolic surface to its length.

\begin{lemma*}[Collar lemma] \label{lem:collar}

If $ \alpha $ is a simple closed geodesic on a hyperbolic surface $ S $,
\[
\len_{ S } \of{ \beta } > i \of{ \alpha , \beta } \cdot 2 \sinh^{ -1 } \of{ \csch \of{ \frac{ \len_{ S } \of{ \alpha } }{ 2 } } }
\]
for any closed curve $ \beta $ on $ S $.

\end{lemma*}

\noindent See \cite[Corollary~4.1.1]{Bus10} for a proof of \autoref{lem:collar}. When all the elevations of a (not necessarily simple) closed geodesic along a finite-degree cover are simple, we will use the above result to bound the self-intersection number in terms of its length.

\begin{corollary*} \label{cor:lengthToIntersection}

Let $ \pi \colon Z \to S $ be a finite-degree cover of a hyperbolic surface $ S $. If every elevation of some curve $ \gamma $ on $ S $ along $ \pi $ to $ Z $ is simple, then its self intersection number $ i \of{ \gamma , \gamma } $ is explicitly bounded in terms of its length $ \len_{ S } \of{ \gamma } $ and the degree $ \deg \of{ \pi } $.

\begin{proof}

Fix a parametrization of $ \gamma $ without triple (or higher order) intersections, and let $ \widetilde{ \gamma }_{ 1 } , \dots , \widetilde{ \gamma }_{ n } $ be the elevations of $ \gamma $ along $ \pi $ to $ Z $, all of which are simple by assumption and have infimal length
\[
\len_{ Z } \of{ \widetilde{ \gamma }_{ k } } \leq \deg \of{ \pi } \cdot \len_{ S } \of{ \gamma } .
\]
Note that $ n \leq \deg \of{ \pi } $. Moreover, every intersection point between elevations of $ \gamma $ lies above a self-intersection point of $ \gamma $, and conversely, the fiber of $ \pi $ over any self-intersection point of $ \gamma $ consists of exactly $ \deg \of{ \pi } $ intersection points between elevations. Hence
\begin{align*}
i \of{ \gamma , \gamma } & = \frac{ 1 }{ 2 \deg \of{ \pi } } \sum_{ j , k = 1 }^{ n }{ i \of{ \widetilde{ \gamma }_{ j } , \widetilde{ \gamma }_{ k } } } < \frac{ 1 }{ 2 \deg \of{ \pi } } \sum_{ j , k = 1 }^{ n }{ \frac{ \len_{ Z } \of{ \widetilde{ \gamma_{ j } } } }{ 2 \sinh^{ -1 } \of{ \csch \of{ \frac{ \len_{ Z } \of{ \widetilde{ \gamma_{ k } } } }{ 2 } } } } } \\
& \leq \frac{ 1 }{ 2 \deg \of{ \pi } } \sum_{ j , k = 1 }^{ n }{ \frac{ \deg \of{ \pi } \len_{ S } \of{ \gamma } }{ 2 \sinh^{ -1 } \of{ \csch \of{ \frac{ \deg \of{ \pi } \len_{ S } \of{ \gamma } }{ 2 } } } } } = \frac{ n^{ 2 } \len_{ S } \of{ \gamma } }{ 4 \sinh^{ -1 } \of{ \csch \of{ \frac{ \deg \of{ \pi } \len_{ S } \of{ \gamma } }{ 2 } } } } \\
& \leq \frac{ \deg \of{ \pi }^{ 2 } \len_{ S } \of{ \gamma } }{ 4 \sinh^{ -1 } \of{ \csch \of{ \frac{ \deg \of{ \pi } \len_{ S } \of{ \gamma } }{ 2 } } } } .
\end{align*}
That is, the self-intersection number $ i \of{ \gamma , \gamma } $ of $ \gamma $ is explicitly bounded in terms of its length $ \len_{ S } \of{ \gamma } $ and the degree $ \deg \of{ \pi } $. \qedhere

\end{proof}

\end{corollary*}

\subsection{Subsurface projection and the bounded geodesic image theorem} \label{sec:subsurfaceProjection}

For this section, fix a basepoint $ s \in S $ on a closed surface $ S $ of negative Euler characteristic. Any curve $ \gamma $ on $ S $ has an associated \emph{annular complex} $ \CC \of{ \gamma } $, which we describe below. 

A curve $ \gamma $ on $ S $ determines a conjugacy class of cyclic subgroups of $ \pi_{ 1 } \of{ S , s } $, and so the Galois correspondence yields a based cover $ \Pi_{ \gamma } \colon \of{ S_{ \gamma } , s_{ \gamma } } \to \of{ S , s } $. The induced homomorphism $ \of{ \Pi_{ \gamma } }_{ * } \colon \pi_{ 1 } \of{ S_{ \gamma } , s_{ \gamma } } \into \pi_{ 1 } \of{ S , s } $ is injective with image the cyclic subgroup $ \gen{ \gamma } $ generated by $ \gamma $, and so $ S_{ \gamma } $ is a surface of finite type with infinite cyclic fundamental group $ \pi_{ 1 } \of{ S_{ \gamma } , s_{ \gamma } } \cong \gen{ \gamma } $. In particular, $ S_{ \gamma } $ is homeomorphic to an open annulus.

Fix a hyperbolic metric on $ S $, which pulls back along $ \Pi_{ \gamma } $ to a complete hyperbolic metric on $ S_{ \gamma } $, and denote by $ \overline{ S_{ \gamma } } $ the Gromov compactification of $ S_{ \gamma } $, which is homeomorphic to a closed annulus. $ \CC \of{ \gamma } $ is a flag simplicial complex whose vertices are path-homotopy classes of properly embedded, non-peripheral arcs in $ \overline{ S_{ \gamma } } $, any pairwise disjoint finite collection of which spans a simplex. As in the case of the curve complex of a closed surface, we equip $ \CC \of{ \gamma } $ with a metric $ d_{ \gamma } $ so that every simplex is a regular Euclidean simplex of side length $ 1 $.

Let $ \mB_{ 2 } \of{ \gamma } $ denote the open ball of radius $ 2 $ in $ \CC \of{ S } $, so that $ \CC \of{ S } \setminus \mB_{ 2 } \of{ \gamma } $ is the sub-complex of $ \CC \of{ S } $ consisting of simplices all of whose vertices intersect $ \gamma $. The pre-image in the cover $S_{\gamma}$ of the geodesic representative of such a simplex $ \alpha \in \CC \of{ S } \setminus \mB_{ 2 } \of{ \gamma } $ contains exactly $ i \of{ \alpha , \gamma } $ properly embedded, non-peripheral geodesic arcs, all of which are disjoint. By choosing one of these arcs, we obtain a coarsely well-defined map $ \proj_{ \gamma } \colon \CC \of{ S } \setminus \mB_{ 2 } \of{ \gamma } \to \CC \of{ \gamma } $, called the \emph{subsurface projection} associated to $ \gamma $.

The following result of Minsky--Taylor gives an explicit bound on the diameter of the subsurface projection onto the annular complex of a non-simple curve.

\begin{theorem*}[{\cite[Theorem~1.3]{MT22}}] \label{thm:subsurfaceProjection}

$ \diam_{ \CC \of{ \gamma } } \of{ \proj_{ \gamma } \of{ \CC \of{ S } \setminus \mB_{ 2 } \of{ \gamma } } } \leq 38 $ for any curve $ \gamma $ on a closed surface $ S $ of negative Euler characteristic which is not a power of a simple curve.

\end{theorem*}

The following corollary follows immediately from the fact that the quasi-isometric embedding of curve complexes induced by a finite-degree cover is distance non-increasing.

\begin{corollary*} \label{cor:subsurfaceProjection}

Let $ p \colon X \to S $ be a finite-degree cover of a closed surface $ S $ of negative Euler characteristic. Then
\[
\diam_{ \CC \of{ \widetilde{ \alpha } } } \of!\Big!{ \proj_{ \widetilde{ \alpha } } \of!\big!{ \widetilde{ p } \of{ \CC \of{ S } } \setminus \mB_{ 2 } \of{ \widetilde{ \alpha } } } } \leq 38
\]
for all simple curves $ \widetilde{ \alpha } \in \CC \of{ X } $ for which $ \alpha \coloneqq p \of{ \widetilde{ \alpha } } $ is not a power of a simple curve.

\begin{proof}

Fix basepoints $ x \in X $ and $ s \in S $ so that $ p \of{ x } = s $. There are choices of cyclic subgroups of $ \pi_{ 1 } \of{ X , x } $ and $ \pi_{ 1 } \of{ S , s } $ representing $ \widetilde{ \alpha } $ and $ \alpha $, respectively, so that the induced homomorphism $ p_{ * } \colon \pi_{ 1 } \of{ X , x } \into \pi_{ 1 } \of{ S , s } $ restricts to an isomorphism of these cyclic subgroups. This implies that there are choices of the corresponding covers $ \Pi_{ \widetilde{ \alpha } } \colon X_{ \widetilde{ \alpha } } \to X $ and $ \Pi_{ \alpha } \colon S_{ \alpha } \to S $ so that $ X_{ \widetilde{ \alpha } } = S_{ \alpha } $ and the following diagram commutes.
\[
\begin{tikzcd}
& X_{ \widetilde{ \alpha } } = S_{ \alpha } \ar[ dl , " \Pi_{ \widetilde \alpha } " swap ] \ar[ dr , " \Pi_{ \alpha } " ] \\
X \ar[ rr , " p " swap ] && S
\end{tikzcd}
\]
Recall that $ \widetilde{ p } $ depends on a choice of elevation of every simple curve on $ S $. If we choose the elevations of curves intersecting $ \alpha $ to intersect $ \widetilde{ \alpha } $, then $ \widetilde{ p } $ restricts to a map $ \CC \of{ S } \setminus \mB_{ 2 } \of{ \alpha } \to \CC \of{ X } \setminus \mB_{ 2 } \of{ \widetilde{ \alpha } } $ so that the following diagram commutes.
\[
\begin{tikzcd}
\CC \of{ S } \setminus \mB_{ 2 } \of{ \alpha } \ar[ dr , " \proj_{ \alpha } " swap ] \ar[ rr , " \widetilde{ p } " ] &&  \CC \of{ X } \setminus \mB_{ 2 } \of{ \widetilde{ \alpha } } \ar[ dl , " \proj_{ \widetilde{ \alpha } } " ] \\
& \CC \of{ \alpha } = \CC \of{ \widetilde{ \alpha } }
\end{tikzcd}
\]
In particular, \autoref{thm:subsurfaceProjection} implies that
\begin{align*}
\diam_{ \CC \of{ \widetilde{ \alpha } } } \of{ \proj_{ \widetilde{ \alpha } } \of{ \widetilde{ p } \of{ \CC \of{ S } } \setminus \mB_{ 2 } \of{ \widetilde{ \alpha } } } } & = \diam_{ \CC \of{ \alpha } } \of{ \proj_{ \alpha } \of{ \CC \of{ S } \setminus \mB_{ 2 } \of{ \alpha } } } \leq 38 . \qedhere
\end{align*}

\end{proof}

\end{corollary*}

We will also make use of the Bounded Geodesic Image theorem, which describes the subsurface projection of a geodesic in the curve complex. Originally proved by Masur--Minsky \cite[Theorem~3.1]{MM00}, we will use the following strengthening due to Webb.

\begin{theorem*}[{\cite[Corollary~1.3]{Web15}}] \label{thm:boundedGeodesicImage}

There is an explicit positive real number $ C > 0 $ so that
\[
\diam_{ \CC \of{ \gamma } } \of{ \proj_{ \gamma } \of{ R } } \leq C
\]
for any simple curve $ \gamma $ on a surface $ S $ of negative Euler characteristic and any geodesic ray $ R $ in $ \CC \of{ S } $ so that $ R \subseteq \mB \of{ \gamma } $.
\end{theorem*}

\noindent It follows from Webb's argument that $ C $ can be taken to be at most $ 1000 $.

\subsection{An explicit computation of \textit{M}}
\label{sub:explicit-computation}

We are now ready to prove \autoref{thm:effective} by computing the integer $ M $ explicitly as a function of the Euler characteristic $ \abs{ \chi \of{ S } } $ and the degrees $ \deg \of{ p } $ and $ \deg \of{ q } $.

Before proceeding, we give a brief sketch of the proof for the reader's convenience. Without loss of generality, throughout we will assume that $ \deg \of{ p } \leq \deg \of{ q } $. First, in \autoref{clm:isomorphicCovers} and \autoref{clm:ifLargerDegree}, we reduce to the case where $ \deg \of{ p } $ and $ \deg \of{ q } $ are equal by using well-known asymptotics of the growth of simple closed curves with respect to length. In particular, we show that if $ \deg \of{ p } < \deg \of{ q } $, then there is always a simple closed curve on $ Y $ whose image in $ S $ has no simple elevation to $ X $, and that the self-intersection number of such a curve can be explicitly bounded. This comprises the simplest part of the argument, and the primary difficulty is confirming that all of the relevant constants are explicitly computable.

In the case that $ \deg \of{ p } = \deg \of{ q } $, for each $ g \in G = \Deck \of{ \pi_{ S } } $, we form an intermediate cover, $ W_{ g } $, corresponding to the subgroup $ H_{ g } $ generated by $ g \Deck \of{ \pi_{ X } } g^{ -1 } $ and $\Deck \of{ \pi_{ Y } } $ (see \autoref{fig:effectiveSunadaDiamond} and the surrounding text). Under the assumption that $ p $ and $ q $ are non-isomorphic covers, it follows that $ Y $ covers each $ W_{ g } $ non-trivially. In \autoref{clm:boundedDistanceCurves}, we then find a simple curve $ \alpha $ on $ X $ and a simple curve $ \eta $ on $ S $ for which we can explicitly compute bounds on the intersection number of elevations $ \widetilde{ \alpha } $ and $ \widetilde{ \eta } $ of $ \alpha $ and $ \eta $ to the cover $ Z $, and for which $ p \of{ \alpha } $ has no simple elevations to $ W_{ g } $ for any $ g \in G $.

Using work of Aougab-Taylor, we then construct a quasi-geodesic ray $ \widetilde{ R } = \of{ \widetilde{ v_{ n } } }_{ n = 0 }^{ \infty } $ in $\CC \of{ Z } $ emanating from $ \widetilde{ \eta } $ where $ \widetilde{ v_{ i } } $ and $ \widetilde{ v_{ j } } $ have explicitly bounded intersection number for each $ i , j \in \N $. The content of \autoref{clm:distanceAtLeast} is then to find some computable universal constant $ E \in \N $ so that twisting $ \widetilde{ R } $ $ E $ times around $ \widetilde{ \alpha } $ makes positive progress away from the curve complexes $\CC \of{ W_{ g } } $; that is, the distance from $ \tau_{ \widetilde{ \alpha } }^{ E } \of{ \widetilde{ v_{ n } } } $ to the quasi-isometrically embedded image of each $ \CC \of{ W_{ g } } $ in $ \CC \of{ Z } $ is bounded below linearly in $ n $. 

Combined with work of Tang, we then use this to show in \autoref{clm:contradiction} that if there is some $ g \in G $ for which $ g \cdot \tau_{ \widetilde{ \alpha } }^{ E } \of{ \widetilde{ v_{ n } } } $ lies in the quasi-isometrically embedded image of $ \CC \of{ Y } $ in $ \CC \of{ Z } $, then the index $ n $ is explicitly bounded. In particular, this claim implies that if $ n $ is large enough, $ g \cdot \tau_{ \widetilde{ \alpha } }^{ E } \of{ \widetilde{ v_{ n } } } $ cannot lie in this image and so for any such $ n $, $ \gamma_{ n } = \pi_{ S } \of{ \tau_{ \widetilde{ \alpha } }^{ E } \of{ \widetilde{ v_{ n } } } } $ will be a curve on $ S $ with a simple elevation to $ X $ but no simple elevation to $ Y $.
Finally in \autoref{clm:finalBound}, we use the work of Aougab-Taylor and Webb to effectively compute an upper bound for the intersection number of one such $ \gamma_{ n } $, completing the proof of \autoref{thm:effective}.

As in the proof of \autoref{cor:awfulCombinatorics}, we emphasize the importance of the existence of certain explicitly computable numbers over that of their exact form computed below. As such, we encourage the reader to skim the proof of each claim on a first reading.

\begin{proof}[Proof of \autoref{thm:effective}]

For ease of reference, let $ \chi \coloneqq \chi \of{ S } $ and $ d \coloneqq \max \of{ \deg \of{ p } , \deg \of{ q } } $. Fix basepoints $ x \in X $, $ y \in Y $, and $ s \in S $ so that $ p \of{ x } = q \of{ y } = s $, and let
\[
K \coloneqq \bigcap_{ \gamma \in \pi_{ 1 } \of{ S , s } }{ \gamma \of{ p_{ * } \of{ \pi_{ 1 } \of{ X , x } } \cap q_{ * } \of{ \pi_{ 1 } \of{ Y , y } } } \gamma^{ -1 } },
\]
be the normal core of $ p_{ * } \of{ \pi_{ 1 } \of{ X , x } } \cap q_{ * } \of{ \pi_{ 1 } \of{ Y , y } } $ in $ \pi_{ 1 } \of{ S , s } $. Then $ K $ is a normal subgroup of $ \pi_{ 1 } \of{ S , s } $ of explicitly bounded index
\begin{align*}
\ind{ \pi_{ 1 } \of{ S , s } }{ K } & \leq \ind{ \pi_{ 1 } \of{ S , s } }{ p_{ * } \of{ \pi_{ 1 } \of{ X , x } } \cap q_{ * } \of{ \pi_{ 1 } \of{ Y , y } } } ! \\
& \leq \of{ \ind{ \pi_{ 1 } \of{ S , s } }{ p_{ * } \of{ \pi_{ 1 } \of{ X , x } } } \cdot \ind{ \pi_{ 1 } \of{ S , s } }{ q_{ * } \of{ \pi_{ 1 } \of{ Y , y } } } } ! \\
& = \of{ \deg \of{ p } \cdot \deg \of{ q } } ! \leq d^{ 2 } !
\end{align*}
which is contained in both $ p_{ * } \of{ \pi_{ 1 } \of{ X , x } } $ and $ q_{ * } \of{ \pi_{ 1 } \of{ Y , y } } $. As in the proof of \cite[Theorem~3]{ALLM23}, there is a surjective homomorphism $ \rho \colon \pi_{ 1 } \of{ S , s } \onto G $ onto a finite group $ G $ with kernel $ \ker \of{ \rho } = K $. Let $ A \coloneqq \rho \of{ p_{ * } \of{ \pi_{ 1 } \of{ X , x } } } $ and $ B \coloneqq \rho \of{ q_{ * } \of{ \pi_{ 1 } \of{ Y , y } } } $ be the corresponding subgroups of $ G $.

Up to re-labeling, we may assume without loss of generality that $ \card{ A } \geq \card{ B } $. For each element $ g \in G $, let $ H_{ g } \coloneqq \gen{ g A g^{ -1 } , B } $. The Galois correspondence between finite-index subgroups of $ \pi_{ 1 } \of{ S , s } $ and finite-degree based covers of $ \of{ S , s } $ then yields for each element $ g \in G $ a based surface $ \of{ W_{ g } , w_{ g } } $ and a commutative diagram of based covers of the form found in \autoref{fig:effectiveSunadaDiamond}.

\begin{figure}[t]
\centering
\begin{tikzcd}
& \pi_{ 1 } \of{ S , s } \ar[ ddl , bend right , no head ] \ar[ d , no head ] \ar[ ddr , bend left , no head ] &&& \of{ Z , z } \ar[ ddd , bend right=50 , " \pi_{ S } " swap ] \ar[ dl , " \pi_{ X } " swap ] \ar[ dd , " \pi_{ g } " swap ] \ar[ dr , " \pi_{ Y } " ] \\
& \rho^{ -1 } \of{ H_{ g } } \ar[ dd , no head ] \ar[ dr , no head ] && \of{ X , x } \ar[ ddr , " p " swap , bend right ] && \of{ Y , y } \ar[ dl , " q_{ g } " ] \ar[ ddl , " q " , bend left ] \\
p_{ * } \of{ \pi_{ 1 } \of{ X , x } } \ar[ dr , no head ] && q_{ * } \of{ \pi_{ 1 } \of{ Y , y } } \ar[ dl , no head ] && \of{ W_{ g } , w_{ g } } \ar[ d , " r_{ g } " swap ] \\
& K &&& \of{ S , s }
\end{tikzcd}
\caption{Covers from the Galois correspondence} \label{fig:effectiveSunadaDiamond}
\end{figure}
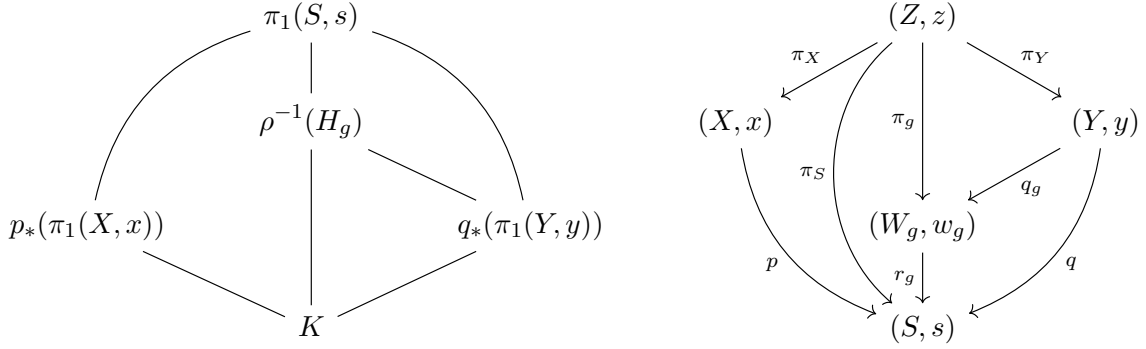

Note that the regular cover $\pi_S$ has explicitly bounded degree
\[
\deg \of{ \pi_{ S } } = \ind{ \pi_{ 1 } \of{ S , s } }{ K } \leq d^{ 2 } !
\]
Moreover, since $ K $ is normal in $ \pi_{ 1 } \of{ S , s } $, this cover is regular with deck group
\[
\Deck \of{ \pi_{ S } } \cong \frac{ N_{ \pi_{ 1 } \of{ S , s } } \of{ K } }{ K } = \frac{ \pi_{ 1 } \of{ S , s } }{ K } \cong \im \of{ \rho } = G .
\]
Similarly, the covers $ \pi_{ X } \colon Z \to X $, $ \pi_{ Y } \colon Z \to Y $, and $ \pi_{ g } \colon Z \to W_{ g } $ are regular, with deck groups $ \Deck \of{ \pi_{ X } } \cong A $, $ \Deck \of{ \pi_{ Y } } \cong B $, $ \Deck \of{ \pi_{ g } } \cong H_{ g } $, respectively. By \cite[Proposition~3.5]{ALLM23}, we may identify the deck groups of these covers by their isomorphic images in $ \Mod \of{ Z } $.

For each element $ g \in G $, the covering relations described in \autoref{thm:CCembedding} yield from the above commutative diagram of finite-degree covering spaces a commutative diagram of coarsely well-defined quasi-isometric embeddings of corresponding curve complexes of the form found in \autoref{fig:quasi-IsometricEmbeddings}.

\begin{figure}[ht]
\centering
\begin{tikzcd}
& \CC \of{ Z } \\
\CC \of{ X } \ar[ ur , " \widetilde{ \pi_{ X } } " ] && \CC \of{ Y } \ar[ ul , " \widetilde{ \pi_{ Y } } " swap ] \\
& \CC \of{ W_{ g } } \ar[ uu , " \widetilde{ \pi_{ g } } " ] \ar[ ur , " \widetilde{ q_{ g } } " swap ] \\
& \CC \of{ S } \ar[ uuu , bend left=50 , " \widetilde{ \pi_{ S } } " ] \ar[ uul , " \widetilde{ p } " , bend left ] \ar[ u , " \widetilde{ r_{ g } } " ] \ar[ uur , " \widetilde{ q } " swap , bend right ]
\end{tikzcd}
\caption{Quasi-isometric embeddings induced by finite-degree covering maps} \label{fig:quasi-IsometricEmbeddings}
\end{figure}
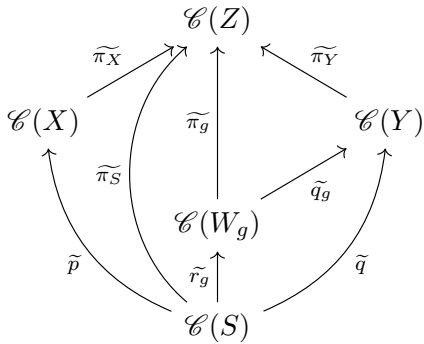

We will later use the fact that an explicitly computable bound on the degree of each cover yields explicit bounds on the quality of the above quasi-isometric embeddings.

\begin{claim} \label{clm:isomorphicCovers}

$ p $ and $ q $ are isomorphic covers if and only if $ B = H_{ g } $ for some element $ g \in G $.

\begin{proof}[Proof of \autoref{clm:isomorphicCovers}]

The Galois correspondence implies that $ p $ and $ q $ are isomorphic as covers of $ S $ precisely when the subgroups $ p_{ * } \of{ \pi_{ 1 } \of{ X , x } } $ and $ q_{ * } \of{ \pi_{ 1 } \of{ Y , y } } $ are conjugate in $ \pi_{ 1 } \of{ S , s } $. Since $G=\Deck(\pi_S)$, this occurs precisely when the images $ \rho \of{ p_{ * } \of{ \pi_{ 1 } \of{ X , x } } } = A $ and $ \rho \of{ q_{ * } \of{ \pi_{ 1 } \of{ Y , y } } } = B $ are conjugate in $ G $. Note that if $ g A g^{ -1 } = B $ for some element $ g \in G $, then
\[
H_{ g } = \gen{ g A g^{ -1 } , B } = \gen{ B , B } = B .
\]
Conversely, if $ H_{ g } = B $ for some element $ g \in G $, then $ g A g^{ -1 } $ is a subgroup of $ H_{ g } = B $ of size $ \card{ g A g^{ -1 } } = \card{ A } \geq \card{ B } $, and so $ g A g^{ -1 } = B $. \qedhere

\end{proof}

\end{claim}

Thus we may assume that $ B \subsetneq H_{ g } $ for every element $ g \in G $, so that $ p $ and $ q $ are not isomorphic as covers of $ S $. In particular, this implies that $ \deg \of{ q_{ g } } = \ind{ H_{ g } }{ B } > 1 $ for each element $ g \in G $. It remains to find an explicitly computable positive integer $ M $ and a curve $ \gamma $ on $ S $ with $ i \of{ \gamma , \gamma } \leq M $ which has a simple elevation to one of $ X $ and $ Y $ but not to both.

To that end, fix an arbitrary hyperbolic metric on $ S $, which pulls back along the covers described above to a hyperbolic metric on each surface in question.

\begin{claim} \label{clm:ifLargerDegree}

If $ \card{ A } > \card{ B } $, then there is an explicitly computable positive integer $ M_{ 1 } $ satisfying the conclusions of \autoref{thm:effective}; that is, $ i \of{ \gamma , \gamma } \leq M_{ 1 } $ for some curve $ \gamma $ on $ S $ which has a simple elevation along $ q $ to $ Y $ but no simple elevations along $ p $ to $ X $.

\begin{proof}[Proof of \autoref{clm:ifLargerDegree}]

Under these hypotheses, note that
\begin{align*}
\deg \of{ q } & = \ind{ \pi_{ 1 } \of{ S } }{ q_{ * } \of{ \pi_{ 1 } \of{ Y , y } } } = \ind{ G }{ B } > \ind{ G }{ A } = \ind{ \pi_{ 1 } \of{ S } }{ p_{ * } \of{ \pi_{ 1 } \of{ X , x } } } = \deg \of{ p } ,
\end{align*}
and so \autoref{cor:awfulCombinatorics} implies that there is an explicitly computable positive number $ L_{ 1 } $ and a simple closed geodesic $ \beta $ on $ Y $ of hyperbolic length at most $ L_{ 1 } $ so that $ \gamma \coloneqq q \of{ \beta } $ has no simple elevations along $ p $ to $ X $. In particular,
\[
L_{ 1 } \coloneqq \max \of{ L_{ 0 } \of{ S } , L_{ 0 } \of{ X } , L_{ 0 } \of{ Y } , 1 + \frac{ c_{ 2 } \of{ X } d^{ 1 + 3 d \abs{ \chi } } }{ c_{ 1 } \of{ Y } } } .
\]
It suffices to find an explicitly computable bound on the self-intersection number $ i \of{ \gamma , \gamma } $. To that end, note that while $ \gamma $ has no simple elevations along $ p $ to $ X $, all of its elevations along $ \pi_{ S } $ to $ Z $ are simple, since $ \beta $ is simple on $ Y $ and $ \pi_{ S } $ is regular. Thus \autoref{cor:lengthToIntersection} implies that $ i \of{ \gamma , \gamma } $ is explicitly bounded; specifically,
\[
i \of{ \gamma , \gamma } \leq \frac{ \deg \of{ \pi_{ S } }^{ 2 } \len_{ S } \of{ \gamma } }{ 4 \sinh^{ -1 } \of{ \csch \of{ \frac{ \deg \of{ \pi_{ S } } \len_{ S } \of{ \gamma } }{ 2 } } } } \leq \frac{ \of{ d^{ 2 } ! }^{ 2 } L_{ 1 } }{ 4 \sinh^{ -1 } \of{ \csch \of{ \frac{ d^{ 2 } ! L_{ 1 } }{ 2 } } } } \eqqcolon M_{ 1 },
\]
as required.
\end{proof}

\end{claim}

Thus we may assume for the remainder of this proof that $ \card{ A } = \card{ B } < \card{ H_{ g } } $ for each element $ g \in G $. Under these assumptions, we will use the full power of \autoref{cor:awfulCombinatorics} to show that $ \widetilde{ \pi_{ X } } \of{ \CC \of{ X } } $ is not contained in $ \bigcup_{ g \in G }{ \widetilde{ \pi_{ g } } \of{ \CC \of{ W_{ g } } } } $, and in fact that the relative complement contains a curve with explicitly bounded length.

\begin{claim} \label{clm:boundedDistanceCurves}

There are simple curves $ \alpha \in \CC \of{ X } $ and $ \eta \in \CC \of{ S } $ so that
\begin{enumerate*}

\item[(i)]

the intersection number $ i \of{ \widetilde{ \pi_{ X } } \of{ \alpha } , \widetilde{ \pi_{ S } } \of{ \eta } } $ is explicitly bounded; and

\item[(ii)]

for all elements $ g \in G $, $ p \of{ \alpha } $ has no simple elevations along $ r_{ g } $ to $ W_{ g } $.

\end{enumerate*}

\begin{proof}[Proof of \autoref{clm:boundedDistanceCurves}]

Since $ \card{ A } < \card{ H_{ g } } $ for each element $ g \in G $,
\begin{align*}
\deg \of{ p } & = \ind{ \pi_{ 1 } \of{ S , s } }{ p_{ * } \of{ \pi_{ 1 } \of{ X , x } } } = \ind{ G }{ A } \\
& > \ind{ G }{ H_{ g } } = \ind{ \pi_{ 1 } \of{ S , s } }{ \of{ r_{ g } }_{ * } \of{ \pi_{ 1 } \of{ W_{ g } , w_{ g } } } } = \deg \of{ r_{ g } } ,
\end{align*}
and so \autoref{cor:awfulCombinatorics} implies that there is a simple closed geodesic $ \alpha $ on $ X $ so that for all elements $ g \in G $, $ p \of{ \alpha } $ has no simple elevations along $ r_{ g } $ to $ W_{ g } $. Moreover, $ \alpha $ can be chosen to be of explicitly bounded length $ \len_{ X } \of{ \alpha } \leq L_{ 2 } $; specifically,
\begin{align*}
\len_{ X } \of{ \alpha } & \leq \max \of{ \max_{ g \in G }{ L_{ 0 } \of{ W_{ g } } } , L_{ 0 } \of{ X } , 1 + \max_{ g \in G }{ \frac{ \card{ G } c_{ 2 } \of{ W_{ g } } \deg \of{ p }^{ 1 + 3 \deg \of{ p } \abs{ \chi } } }{ c_{ 1 } \of{ X } } } } \\
& \leq \max \of{ \max_{ g \in G }{ L_{ 0 } \of{ W_{ g } } } , L_{ 0 } \of{ X } , 1 + \frac{ d^{ 2 } ! \of{ \max_{ g \in G }{ c_{ 2 } \of{ W_{ g } } } } d^{ 1 + 3 d \abs{ \chi } } }{ c_{ 1 } \of{ X } } } \eqqcolon L_{ 2 } .
\end{align*}
Explicit bounds on $ \deg \of{ \pi_{ X } } $ and $ \len_{ X } \of{ \alpha } $ now give explicit bounds on the elevation $ \widetilde{ \pi_{ X } } \of{ \alpha } $; specifically,
\begin{equation} \label{eqn:gnarlyal}
    \len_{ Z } \of{ \widetilde{ \pi_{ X } } \of{ \alpha } } \leq \deg \of{ \pi_{ X } } \cdot \len_{ X } \of{ \alpha } \leq \deg \of{ \pi_{ S } } \cdot \len_{ X } \of{ \alpha } \leq d^{ 2 } ! \cdot L_{ 2 }.
\end{equation}
By Bers' theorem (see \cite[Theorem~5.1.2]{Bus10}), there is a simple closed geodesic $ \eta \in S $ of length at most $ 4 \log \of{ 4 \pi \abs{ \chi } } $. A similar argument now gives an explicit bound on the length of the elevation $ \widetilde{ \pi_{ S } } \of{ \eta } $. Explicitly,
\[
\len_{ Z } \of{ \widetilde{ \pi_{ S } } \of{ \eta } } \leq \deg \of{ \pi_{ S } } \cdot \len_{ S } \of{ \eta } \leq d^{ 2 } ! \cdot 4 \log \of{ 4 \pi \abs{ \chi } } .
\]
\autoref{lem:collar} now gives the following explicit bound on the intersection number $ i \of{ \widetilde{ \pi_{ X } } \of{ \alpha } , \widetilde{ \pi_{ S } } \of{ \eta } } $:
\[
i \of{ \widetilde{ \pi_{ X } } \of{ \alpha } , \widetilde{ \pi_{ S } } \of{ \eta } } \leq \frac{ d^{ 2 } ! \cdot 4 \log \of{ 4 \abs{ \chi } } }{ 2 \sinh^{ -1 } \of{ \csch \of{ \frac{ d^{ 2 } ! L_{ 2 } }{ 2 } } } } \coloneqq K_{ 1 } , 
\]
as required.
\end{proof}

\end{claim}

Fix curves $ \alpha $ and $ \eta $ satisfying the conclusions of \autoref{clm:boundedDistanceCurves}. For ease of notation, we will denote by $ \widetilde{ \alpha } $ and $ \widetilde{ \eta } $ the elevations $ \widetilde{ \pi_{ X } } \of{ \alpha } $ and $ \widetilde{ \pi_{ S } } \of{ \eta } $, respectively. Note that Hempel's bound of \cite[Lemma~2.1]{Hem01} gives an explicit bound
\begin{equation}\label{eqn:logcity}
d_{ Z } \of{ \widetilde{ \alpha } , \widetilde{ \eta } } \leq 2 \log_{ 2 } \of{ i \of{ \widetilde{ \alpha } , \widetilde{ \eta } } } + 2 \leq 2 \log_{ 2 } \of{ K_{ 1 } } + 2,
\end{equation}
on the distance $ d_{ Z } \of{ \widetilde{ \alpha } , \widetilde{ \eta } } $. Moreover, for each element $ g \in G $, since $ \pi_{ g } \of{ \widetilde{ \alpha } } $ is an elevation of $ p \of{ \alpha } $ along $ r_{ g } $ to $ W_{ g } $, it is non-simple. \autoref{cor:subsurfaceProjection} therefore implies that
$$ \diam_{ \CC \of{ \widetilde{ \alpha } } } \of{ \proj_{ \widetilde{ \alpha } } \of{ \widetilde{ \pi_{ g } } \of{ \CC \of{ W_{ g } } } \setminus \mB_{ 2 } \of{ \widetilde{ \alpha } } } } \leq 38. $$
Let $ \of{ v_{ k } }_{ k = 0 }^{ \infty } $ be a sequence in $ \CC \of{ S } $ with $ v_{ 0 } = \eta $ and
\begin{equation} \label{eqn:gnarlychi}
   d_{ S } \of{ v_{ j } , v_{ k } }  = \abs{ j - k }, \quad i \of{ v_{ j } , v_{ k } } \leq \of{ C + 3 }^{ 2 \max \of{ j , k } - 5 } \of{ \frac{ 3 \abs{ \chi } }{ 2 } }^{ \abs{ j - k } + 2 } + f_{ j , k } \of{ \abs{ \chi } }, 
\end{equation}
for all indices $ j , k \in \N $, whose existence is guaranteed by \autoref{thm:aougab-taylor}. We recall that the constant $ C $ appearing above is the fixed, computable universal constant appearing in \autoref{thm:boundedGeodesicImage}.

For ease of reference, we will denote
\begin{align*}
R & \coloneqq \of{ v_{ n } }_{ n = 0 }^{ \infty } & \widetilde{ v_{ n } } & \coloneqq \widetilde{ \pi_{ S } } \of{ v_{ n } } & \widetilde{ R } & \coloneqq \of{ \widetilde{ v_{ n } } }_{ n = 0 }^{ \infty } = \of{ \widetilde{ \pi_{ S } } \of{ v_{ n } } }_{ n = 0 }^{ \infty } ,
\end{align*}
and we will denote by $ \tau_{ \widetilde{ \alpha } } \in \Mod \of{ Z } $ a Dehn twist around $ \widetilde{ \alpha } $.

\begin{claim} \label{clm:distanceAtLeast}

There is a computable universal constant $ E \in \N $ and an explicitly computable constant $K_{2} \geq 1$ so that
\[
d_{ Z } \of{ \tau_{ \widetilde{ \alpha } }^{ E } \of{ \widetilde{ v_{ n } } } , \widetilde{ \pi_{ g } } \of{ \CC \of{ W_{ g } } } } \geq  \frac{n}{K_{2}} - K_{2},
\]
for all indices $ n \in \N $ and elements $ g \in G $.

\begin{proof}[Proof of \autoref{clm:distanceAtLeast}]

We will find some $E$ so that for each $n$ and for any $z \in \CC \of{ W_{ g } } $, any geodesic segment $[\tau^{E}_{\widetilde{\alpha}}(\widetilde{v_{n}}), \widetilde{\pi_{g}}(z)]$ passes within distance $1$ of $\widetilde{\alpha}$ in $\CC(Z)$. The claim will follow since the confluence of \autoref{eqn:logcity} and \autoref{eqn:gnarlychi} show that $\widetilde{\alpha}$ is at distance at least $|n- 2\log_{2}(K_{1}) -2|$ from $\tau^{E}_{\widetilde{\alpha}}(\widetilde{v_{n}})$, where the multiplicative constant comes from the quasi-geodesic constants in \autoref{thm:APT22} which depend explicitly on the degrees of the covers. 

Let $ E \coloneqq C+41$ and for ease of notation let $\widetilde{z}=\widetilde{\pi_g}(z)$.
For contradiction, we will assume that a geodesic segment connecting $\tau_{\widetilde{\alpha}}^E(\widetilde{v_n})$ and $\widetilde{z}$ does not pass within distance $1$ of $\widetilde{\alpha}$ from which we will violate the bound in \autoref{thm:boundedGeodesicImage}.
Note that \autoref{cor:subsurfaceProjection} implies that
\[
\diam_{ \CC \of{ \widetilde{ \alpha } } } \of{ \proj_{ \widetilde{ \alpha } } \of{ \widetilde{ \pi_{ g } } \of{ \CC \of{ W_{ g } } } } } = \diam_{ \CC \of{ \widetilde{ \alpha } } } \of{ \proj_{ \widetilde{ \alpha } } \of{ \widetilde{ \pi_{ g } } \of{ \CC \of{ W_{ g } } } \setminus \mB_{ 2 } \of{ \widetilde{ \alpha } } } } \leq 38 .
\] 
Since given any $\kappa \in \CC \of{Z}$ which lies outside of the $1$-neighborhood of $\widetilde{\alpha}$, the distance in the annular complex of $\widetilde{\alpha}$ between $\kappa$ and $\tau^{E}_{\widetilde{\alpha}}(\kappa)$ is at least $E-2$, one has that 
\begin{align*} d_{\CC \of{\widetilde{\alpha}}}(\proj_{\widetilde{\alpha}}(\tau^{E}_{\widetilde{\alpha}}(\widetilde{v_n})), \proj_{\widetilde{\alpha}}(\widetilde{z}))&\geq d_{\CC \of{\widetilde{\alpha}}}(\proj_{\widetilde{\alpha}}(\tau^{E}_{\widetilde{\alpha}}(\widetilde{v_n})), \proj_{\widetilde{\alpha}}(\widetilde{v_n}))-d_{\CC \of{\widetilde{\alpha}}}(\proj_{\widetilde{\alpha}}(\widetilde{v_n}),\proj_{\widetilde{\alpha}}(\widetilde{z}) )\\ &\geq E- 40, \end{align*}
where the final line follows since $\widetilde {v}_n,\widetilde{z}\in\widetilde{\pi_g}(\CC(W_g))$.
By choice of $E$, this lower bound is at least $C+1$, contradicting \autoref{thm:boundedGeodesicImage}.
Therefore any geodesic connecting $\tau^{E}_{\widetilde{\alpha}}(\widetilde{v_{n}})$ to $\widetilde{z}$ passes through the $1$-neighborhood of $\widetilde{\alpha}$, as desired. 
\end{proof}

\end{claim}

For each index $ n \in \N $, let $ \tau_{ n } \coloneqq \tau_{ \widetilde{ \alpha } }^{ E } \of{ \widetilde{ v_{ n } } } $.
Note that $\widetilde{\alpha},\widetilde{v_n}\in \widetilde{\pi_X}(\CC(X))$ and therefore $\tau_n\in\widetilde{\pi_X}(\CC(X))$.
Moreover, for each element $ g \in G $, consider the nearest point projection map $ \nu_{ g } \colon \CC \of{ Z } \onto \widetilde{ \pi_{ g } } \of{ \CC \of{ W_{ g } } } $, which is coarsely well-defined, and let $ \beta_{ g } \of{ \tau_{ n } } $ denote a circumcenter of the $ H_{ g } $-orbit of $ \tau_{ n } $. \autoref{thm:Tang} implies that
\[
d_{ Z } \of{ \beta_{ g } \of{ \tau_{ n } } , \nu_{ g } \of{ \tau_{ n } } } \leq D,
\]
for some positive real number $ D $ which, as we will show in \autoref{sec:effectivizingTang}, can be explicitly bounded in terms of $ \abs{ \chi } $ and $ d $. In particular, note that
 
\begin{align*}
\frac{n}{K_{2}} - K_{2} & \leq d_{ Z } \of{ \tau_{ n } , \widetilde{ \pi_{ g } } \of{ \CC \of{ W_{ g } } } } = d_{ Z } \of{ \tau_{ n } , \nu_{ g } \of{ \tau_{ n } } } \\
& \leq d_{ Z } \of{ \tau_{ n } , \beta_{ g } \of{ \tau_{ n } } } + d_{ Z } \of{ \beta_{ g } \of{ \tau_{ n } } , \nu_{ g } \of{ \tau_{ n } } } \leq \diam_{ \CC \of{ Z } } \of{ H_{ g } \cdot \tau_{ n } } + D ,
\end{align*}
so that $ \diam_{ \CC \of{ Z } } \of{ H_{ g } \cdot \tau_{ n } } \geq \frac{n}{K_{2}} - K_{2} - D $.

\begin{claim} \label{clm:contradiction}

If $ g \cdot \tau_{ n } \in \widetilde{ \pi_{ Y } } \of{ \CC \of{ Y } } $ for some element $ g \in G $, then $ n $ is explicitly bounded; specifically, $ n \leq K_{ 2 } \of{ d^{ 2 } ! + K_{ 2 } + D } $.

\begin{proof}[Proof of \autoref{clm:contradiction}]

Since $ g \cdot \tau_{ n } \in \widetilde{ \pi_{ Y } } \of{ \CC \of{ Y } } $, its orbit under $ B = \Deck \of{ \pi_{ Y } } $ is a multi-curve on $ Z $; that is, $ d_{ Z } \of{ b \cdot \of{ g \cdot \tau_{ n } } , g \cdot \tau_{ n } } \leq 1 $ for all elements $ b \in B $. Similarly, since $ \tau_{ n } \in \widetilde{ \pi_{ X } } \of{ \CC \of{ X } } $, its orbit under $ A = \Deck \of{ \pi_{ X } } $ is a multi-curve on $ Z $, which implies that
\[
d_{ Z } \of{ \of{ g a g^{ -1 } } \cdot \of{ g \cdot \tau_{ n } } , g \cdot \tau_{ n } } = d_{ Z } \of{ g a \cdot \tau_{ n } , g \cdot \tau_{ n } } = d_{ Z } \of{ a \cdot \tau_{ n } , \tau_n } \leq 1 .
\]
Since $ \of{ g A g^{ -1 } } \cup B $ generates $ H_{ g } $, every element $ h \in H_{ g } $ can be written as a word in $ \of{ g A g^{ -1 } } \cup B $ of length at most $ \card{ H_{ g } } $. In particular,
\[
d_{ Z } \of{ h \cdot \of{ g \cdot \tau_{ n } } , g \cdot \tau_{ n } } \leq d_{ \of{ g A g^{ -1 } } \cup B } \of{ 1 , h } \leq \card{ H_{ g } } \leq \card{ G } \leq d^{ 2 } !
\]
for all elements $ h \in H_{ g } $, and so
\[
n \leq K_{ 2 } \cdot \diam_{ \CC \of{ Z } } \of{ H_{ g } \cdot \tau_{ n } } + K_{ 2 }^{ 2 } + D K_{ 2 } \leq K_{ 2 } \cdot \card{ G } + K_{ 2 }^{ 2 } + D K_{ 2 } \leq K_{ 2 } \of{ d^{ 2 } ! + K_{ 2 } + D } ,
\]
as required. \qedhere

\end{proof}

\end{claim}

Fix an explicitly computable natural number $ K_{ 3 } \in \N $ so that $ K_{ 2 } \cdot \of{ d^{ 2 } ! + K_{ 2 } + D } + 1 < 2^{ K_{ 3 } } + 2 $. Let $ N \coloneqq 2^{ K_{ 3 } } + 2 $, and note that \autoref{clm:contradiction} implies that $ \tau_{ N } \in \widetilde{ \pi_{ X } } \of{ \CC \of{ X } } $ but $ g \cdot \tau_{ N } \notin \widetilde{ \pi_{ Y } } \of{ \CC \of{ Y } } $ for all elements $ g \in G $. In particular, $ \gamma \coloneqq \pi_{ S } \of{ \tau_{ N } } $ has a simple elevation to $ X $ but no simple elevations to $ Y $. It therefore suffices to show that the self-intersection number $ i \of{ \gamma , \gamma } $ is explicitly bounded. In fact, we will compute an explicit bound on the length $ \len_{ S } \of{ \gamma } $ and again argue via the collar lemma.

\begin{claim} \label{clm:finalBound}

The self-intersection number $ i \of{ \gamma , \gamma } $ is explicitly bounded; that is, $ i \of{ \gamma , \gamma } \leq M_{ 2 } $ for some explicitly computable positive integer $ M_{ 2 } $.

\begin{proof}[Proof of \autoref{clm:finalBound}]

Continuing to let $(v_k)_{k=0}^\infty$ denote the sequence from Theorem 2.2, we point out that the construction of the $v_k$ is via an iterated Dehn twisting argument.
In particular, following \cite[Section 5]{AT14}, the $v_k$ are defined so that $v_4=T^{C+3}_{v_3}(v_0)$ and for each $k\ge 2$, $v_{2^k+2}=T^{C+3}_{v^{2^{k-1}+2}}(v_0)$ where $C$ is the constant coming from the bounded geodesic image theorem. The constant
$C$ can be taken to be less than $938$ by work of Webb \cite{Web15}.
One can then compute inductively that 
\begin{align*}
\ell_S(v_4)&=\ell_S(T^{941}_{v_3}(v_0))\le941~\!i(v_3,v_0)\ell_S(v_3)+\ell_S(v_0),\\
\ell_S(v_{6})&=\ell_S(T^{941}_{v_4}(v_0))\le941~\!i(v_4,v_0)\ell_S(v_4)+\ell_S(v_0),\\
&\vdots\\
\ell_S(v_{N_K})=\ell_S(v_{2^K+2})&=\ell_S(T^{941}_{v_{2^{K-1}+2}}(v_0))\le941~\!i(v_{2^{K-1}+2},v_0)\ell_S(v_{2^{K-1}+2})+\ell_S(v_0).
\end{align*}
Using the fact that $i(v_k,v_j)$ can be explicitly bounded using \autoref{eqn:gnarlychi} we conclude inductively from the above equation that there is an effectively computable upper bound on $\ell_S(v_{N_K})$ with leading term of the form
\begin{align*}941~\!i(v_3,v_0)\ell_S(v_3)\prod_{k=1}^{K-1}941i(v_{2^k+2},v_0)&\le 941^{K+1}\left(\frac{3|\chi|}{2}\right)^{5}\ell_S(v_3)\prod_{k=1}^{K-1}\left(941^{2^{k+1}-1}\left(\frac{3|\chi|}{2}\right)^{2^k+4}\right),\\
&=941^{2^{K+1}-2}\left(\frac{3|\chi|}{2}\right)^{2^K+4K-1},
\end{align*}
where we assume that the product on the left side of the first equation is empty whenever $K< 2$.
Consequently the same holds of $\ell_{N_K}=\ell_Z(\widetilde{\pi_{S}}(v_{N_K}))$ since $\ell_Z(\widetilde{\pi_{S}}(v_{N_K}))\le |G|~\!\ell_S(v_{N_K})$.
It then follows that
\begin{equation}\label{eqn:lengthbound}
\ell_Z(\tau^E_{\widetilde{\alpha}}(\widetilde{\pi_{S}}(v_{N_K})))\le E~\!i(\widetilde{\pi_{S}}(v_{N_K}),\widetilde{\alpha})\ell_Z(\widetilde\alpha)+\ell_{N_K}.
\end{equation}
We have an effective bound on the length of $\widetilde{\alpha}$ via \autoref{eqn:gnarlyal}.
Using \autoref{lem:collar} it therefore follows that we can effectively compute an upper bound on $i(\widetilde{\pi_{S}}(v_{N_K}),\widetilde\alpha)$ and hence all of the terms on the righthand side of \autoref{eqn:lengthbound} are effectively computable and therefore so is a bound on $\ell_Z(\tau^E_{\widetilde{\alpha}}(\widetilde{\pi_{S}}(v_{N_K})))$.
Since $\gamma=\pi_{S}(\tau^E_{\widetilde{\alpha}}(\widetilde{\pi_{S}}(v_{N_K})))$ it follows that
$$\ell_S(\gamma)\le\ell_Z(\tau^E_{\widetilde{\alpha}}(\widetilde{\pi_{S}}(v_{N_K}))).$$
We are therefore done since the collar lemma again provides an effectively computable upper bound for $i(\gamma,\gamma)$.
\end{proof}

\end{claim}

Finally, let us summarize what we have shown. If $ p $ and $ q $ are not isomorphic as covers of $ S $, then there are two cases; either $ \card{ A } = \card{ B } $, or $ \card{ A } > \card{ B } $. If $ \card{ A } = \card{ B } $, we showed that $ i \of{ \gamma , \gamma } \leq M_{ 1 } $ for some curve $ \gamma $ on $ S $ which has a simple elevation along $ q $ to $ Y $ but no simple elevations along $ p $ to $ X $. On the other hand, if $ \card{ A } > \card{ B } $, we showed that $ i \of{ \gamma , \gamma } \leq M_{ 2 } $ for some curve $ \gamma $ on $ S $ which has a simple elevation along $ p $ to $ X $ but no simple elevations along $ q $ to $ Y $. Therefore, $ M \coloneqq \max \of{ M_{ 1 } , M_{ 2 } } $ satisfies the conclusions of \autoref{thm:effective}. \qedhere 

\setcounter{claim}{0}

\end{proof}

We remark that the above proof of \autoref{thm:effective} shows that so long as $ \alpha \in \CC \of{ X } $ is sufficiently far from the pullback of each $ \CC \of{ W_{ g } } $ in $ \CC \of{ Z } $, it cannot lie over a curve in $ S $ that lifts simply to $ Y $. Indeed, being far away from each $ \CC \of{ W_{ g } } $ guarantees that the $ \Deck \of{ \pi_{ g } } $-orbits of $ \widetilde{ \alpha } $ are of large diameter. On the other hand, we show that some of these orbits are relatively small when $ p \of{ \alpha } $ has simple elevations along $ q $ to $ Y $. Thus, in the case that the covers in question are not isomorphic, our proof of \autoref{thm:effective} exhibits not only one curve on $ S $ with bounded self-intersection and which lifts simply to one cover but not the other, but in fact shows that \emph{``almost all''} curves on $ S $ that lift simply to one cover do not lift simply to the other.

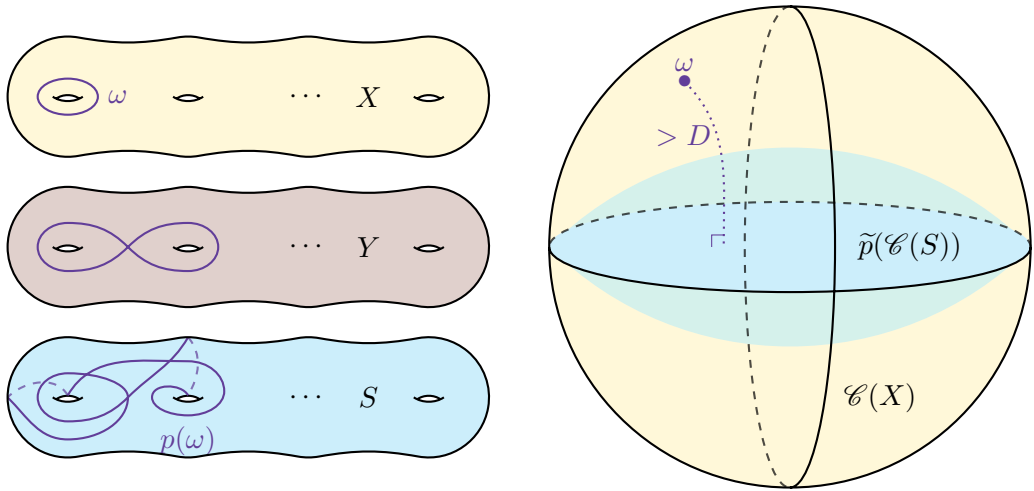
\begin{figure}[ht]
\centering
\begin{tikzpicture}[ scale = 0.8 ]

\fill[ Goldenrod , opacity = 0.2 ]
    ( 0 , 1 + 2.5 ) arc[ radius = 1 , start angle = 90 , end angle = { 270 + 10 } ]
        to[ out = 10 , in = { 180 - 10 } ] ( 2 - 0.1736 , -0.9848 + 2.5 )
        arc[ radius = 1 , start angle = { -90 - 10 } , end angle = { -90 + 10 } ]
        to[ out = 10 , in = { 180 - 10 } ] ( 4 - 0.1736 , -0.9848 + 2.5 )
        arc[ radius = 1 , start angle = { -90 - 10 } , end angle = { -90 + 10 } ]
        to[ out = 10 , in = { 180 - 10 } ] ( 6 - 0.1736 , -0.9848 + 2.5 )
        arc[ radius = 1 , start angle = { -90 - 10 } , end angle = { 90 + 10 } ]
        to[ out = { 180 + 10 } , in = -10 ] ( 4 + 0.1736 , 0.9848 + 2.5 )
        arc[ radius = 1 , start angle = { 90 - 10 } , end angle = { 90 + 10 } ]
        to[ out = { 180 + 10 } , in = -10 ] ( 2 + 0.1736 , 0.9848 + 2.5 )
        arc[ radius = 1 , start angle = { 90 - 10 } , end angle = { 90 + 10 } ]
        to[ out = { 180 + 10 } , in = -10 ] ( 0 + 0.1736 , 0.9848 + 2.5 )
        arc[ radius = 1 , start angle = { 90 - 10 } , end angle = 90 ] ;
\fill[ Sepia , opacity = 0.2 ]
    ( 0 , 1 ) arc[ radius = 1 , start angle = 90 , end angle = { 270 + 10 } ]
        to[ out = 10 , in = { 180 - 10 } ] ( 2 - 0.1736 , -0.9848 )
        arc[ radius = 1 , start angle = { -90 - 10 } , end angle = { -90 + 10 } ]
        to[ out = 10 , in = { 180 - 10 } ] ( 4 - 0.1736 , -0.9848 )
        arc[ radius = 1 , start angle = { -90 - 10 } , end angle = { -90 + 10 } ]
        to[ out = 10 , in = { 180 - 10 } ] ( 6 - 0.1736 , -0.9848 )
        arc[ radius = 1 , start angle = { -90 - 10 } , end angle = { 90 + 10 } ]
        to[ out = { 180 + 10 } , in = -10 ] ( 4 + 0.1736 , 0.9848 )
        arc[ radius = 1 , start angle = { 90 - 10 } , end angle = { 90 + 10 } ]
        to[ out = { 180 + 10 } , in = -10 ] ( 2 + 0.1736 , 0.9848 )
        arc[ radius = 1 , start angle = { 90 - 10 } , end angle = { 90 + 10 } ]
        to[ out = { 180 + 10 } , in = -10 ] ( 0 + 0.1736 , 0.9848 )
        arc[ radius = 1 , start angle = { 90 - 10 } , end angle = 90 ] ;
\fill[ Cyan , opacity = 0.2 ]
    ( 0 , 1 - 2.5 ) arc[ radius = 1 , start angle = 90 , end angle = { 270 + 10 } ]
        to[ out = 10 , in = { 180 - 10 } ] ( 2 - 0.1736 , -0.9848 - 2.5 )
        arc[ radius = 1 , start angle = { -90 - 10 } , end angle = { -90 + 10 } ]
        to[ out = 10 , in = { 180 - 10 } ] ( 4 - 0.1736 , -0.9848 - 2.5 )
        arc[ radius = 1 , start angle = { -90 - 10 } , end angle = { -90 + 10 } ]
        to[ out = 10 , in = { 180 - 10 } ] ( 6 - 0.1736 , -0.9848 - 2.5 )
        arc[ radius = 1 , start angle = { -90 - 10 } , end angle = { 90 + 10 } ]
        to[ out = { 180 + 10 } , in = -10 ] ( 4 + 0.1736 , 0.9848 - 2.5 )
        arc[ radius = 1 , start angle = { 90 - 10 } , end angle = { 90 + 10 } ]
        to[ out = { 180 + 10 } , in = -10 ] ( 2 + 0.1736 , 0.9848 - 2.5 )
        arc[ radius = 1 , start angle = { 90 - 10 } , end angle = { 90 + 10 } ]
        to[ out = { 180 + 10 } , in = -10 ] ( 0 + 0.1736 , 0.9848 - 2.5 )
        arc[ radius = 1 , start angle = { 90 - 10 } , end angle = 90 ] ;

\fill[ White ]
    ( 0 - 0.1926 , -0.0329 + 2.5 ) arc[ x radius = 0.2723 , y radius = 0.2 , start angle = 135 , end angle = 45 ]
    arc[ x radius = 0.35355 , y radius = 0.25 , start angle = 303 , end angle = 237 ] ;
\draw[ thick ]
    ( 0 - 0.25 , 0 + 2.5 ) arc[ x radius = 0.35355 , y radius = 0.25 , start angle = 225 , end angle = 315 ]
    ( 0 - 0.1926 , -0.0329 + 2.5 ) arc[ x radius = 0.2723 , y radius = 0.2 , start angle = 135 , end angle = 45 ] ;

\fill[ White ]
    ( 2 - 0.1926 , -0.0329 + 2.5 ) arc[ x radius = 0.2723 , y radius = 0.2 , start angle = 135 , end angle = 45 ]
    arc[ x radius = 0.35355 , y radius = 0.25 , start angle = 303 , end angle = 237 ] ;
\draw[ thick ]
    ( 2 - 0.25 , 0 + 2.5 ) arc[ x radius = 0.35355 , y radius = 0.25 , start angle = 225 , end angle = 315 ]
    ( 2 - 0.1926 , -0.0329 + 2.5 ) arc[ x radius = 0.2723 , y radius = 0.2 , start angle = 135 , end angle = 45 ] ;

\draw ( 4 , 0 + 2.5 ) node {$ \cdots $} ;

\fill[ White ]
    ( 6 - 0.1926 , -0.0329 + 2.5 ) arc[ x radius = 0.2723 , y radius = 0.2 , start angle = 135 , end angle = 45 ]
    arc[ x radius = 0.35355 , y radius = 0.25 , start angle = 303 , end angle = 237 ] ;
\draw[ thick ]
    ( 6 - 0.25 , 0 + 2.5 ) arc[ x radius = 0.35355 , y radius = 0.25 , start angle = 225 , end angle = 315 ]
    ( 6 - 0.1926 , -0.0329 + 2.5 ) arc[ x radius = 0.2723 , y radius = 0.2 , start angle = 135 , end angle = 45 ] ;

\fill[ White ]
    ( 0 - 0.1926 , -0.0329 ) arc[ x radius = 0.2723 , y radius = 0.2 , start angle = 135 , end angle = 45 ]
    arc[ x radius = 0.35355 , y radius = 0.25 , start angle = 303 , end angle = 237 ] ;
\draw[ thick ]
    ( 0 - 0.25 , 0 ) arc[ x radius = 0.35355 , y radius = 0.25 , start angle = 225 , end angle = 315 ]
    ( 0 - 0.1926 , -0.0329 ) arc[ x radius = 0.2723 , y radius = 0.2 , start angle = 135 , end angle = 45 ] ;

\fill[ White ]
    ( 2 - 0.1926 , -0.0329 ) arc[ x radius = 0.2723 , y radius = 0.2 , start angle = 135 , end angle = 45 ]
    arc[ x radius = 0.35355 , y radius = 0.25 , start angle = 303 , end angle = 237 ] ;
\draw[ thick ]
    ( 2 - 0.25 , 0 ) arc[ x radius = 0.35355 , y radius = 0.25 , start angle = 225 , end angle = 315 ]
    ( 2 - 0.1926 , -0.0329 ) arc[ x radius = 0.2723 , y radius = 0.2 , start angle = 135 , end angle = 45 ] ;

\draw ( 4 , 0 ) node {$ \cdots $} ;

\fill[ White ]
    ( 6 - 0.1926 , -0.0329 ) arc[ x radius = 0.2723 , y radius = 0.2 , start angle = 135 , end angle = 45 ]
    arc[ x radius = 0.35355 , y radius = 0.25 , start angle = 303 , end angle = 237 ] ;
\draw[ thick ]
    ( 6 - 0.25 , 0 ) arc[ x radius = 0.35355 , y radius = 0.25 , start angle = 225 , end angle = 315 ]
    ( 6 - 0.1926 , -0.0329 ) arc[ x radius = 0.2723 , y radius = 0.2 , start angle = 135 , end angle = 45 ] ;

\fill[ White ]
    ( 0 - 0.1926 , -0.0329 - 2.5 ) arc[ x radius = 0.2723 , y radius = 0.2 , start angle = 135 , end angle = 45 ]
    arc[ x radius = 0.35355 , y radius = 0.25 , start angle = 303 , end angle = 237 ] ;
\draw[ thick ]
    ( 0 - 0.25 , 0 - 2.5 ) arc[ x radius = 0.35355 , y radius = 0.25 , start angle = 225 , end angle = 315 ]
    ( 0 - 0.1926 , -0.0329 - 2.5 ) arc[ x radius = 0.2723 , y radius = 0.2 , start angle = 135 , end angle = 45 ] ;

\fill[ White ]
    ( 2 - 0.1926 , -0.0329 - 2.5 ) arc[ x radius = 0.2723 , y radius = 0.2 , start angle = 135 , end angle = 45 ]
    arc[ x radius = 0.35355 , y radius = 0.25 , start angle = 303 , end angle = 237 ] ;
\draw[ thick ]
    ( 2 - 0.25 , 0 - 2.5 ) arc[ x radius = 0.35355 , y radius = 0.25 , start angle = 225 , end angle = 315 ]
    ( 2 - 0.1926 , -0.0329 - 2.5 ) arc[ x radius = 0.2723 , y radius = 0.2 , start angle = 135 , end angle = 45 ] ;

\draw ( 4 , 0 - 2.5 ) node {$ \cdots $} ;

\fill[ White ]
    ( 6 - 0.1926 , -0.0329 - 2.5 ) arc[ x radius = 0.2723 , y radius = 0.2 , start angle = 135 , end angle = 45 ]
    arc[ x radius = 0.35355 , y radius = 0.25 , start angle = 303 , end angle = 237 ] ;
\draw[ thick ]
    ( 6 - 0.25 , 0 - 2.5 ) arc[ x radius = 0.35355 , y radius = 0.25 , start angle = 225 , end angle = 315 ]
    ( 6 - 0.1926 , -0.0329 - 2.5 ) arc[ x radius = 0.2723 , y radius = 0.2 , start angle = 135 , end angle = 45 ] ;

\draw[ thick , RoyalPurple ]
    ( 0 , 2.5 ) circle[ x radius = 0.5 , y radius = 0.3 ]
    ( 0.5 , 2.5 ) node[ anchor = west ] {$ \omega $}
    ( 1 , 0 ) to[ out = 45 , in = 180 ] ( 2 , 0.4 )
        arc[ x radius = 0.5 , y radius = 0.4 , start angle = 90 , end angle = -90 ]
        to[ out = 180 , in = 315 ] ( 1 , 0 )
        to[ out = 135 , in = 0 ] ( 0 , 0.4 )
        arc[ x radius = 0.5 , y radius = 0.4 , start angle = 90 , end angle = 270 ]
        to[ out = 0 , in = 225 ] ( 1 , 0 )
    ( 1 , -2.5 ) arc[ x radius = 1 , y radius = 0.4 , start angle = 0 , end angle = 90 ]
        arc[ x radius = 0.5 , y radius = 0.4 , start angle = 90 , end angle = 270 ]
        to[ out = 0 , in = 225 ] ( 1 , -2.5 )
        to[ out = 45 , in = 240 ] ( 2 , -2.5 + 1 )
    ( 2 , -2.5 + 0.0329 ) to[ out = 120 , in = 90 ] ( 1.4 , -2.5 )
        arc[ x radius = 0.6 , y radius = 0.3 , start angle = 180 , end angle = 270 ]
            node[ anchor = north ] {$ p \of{ \omega } $}
        arc[ x radius = 0.6 , y radius = 0.4 , start angle = 270 , end angle = 360 ]
        arc[ x radius = 0.6 , y radius = 0.5 , start angle = 0 , end angle = 90 ]
        to[ out = 180 , in = 60 ] ( 0 , -2.5 + 0.0329 )
    ( -1 , -2.5 ) to[ out = 315 , in = 180 ] ( 0 , -2.5 - 0.7 )
        arc[ x radius = 1 , y radius = 0.7 , start angle = 270 , end angle = 360 ] ;
\draw[ thick , dashed , RoyalPurple , opacity = 0.7 ]
    ( 2 , -2.5 + 1 ) to[ out = 310 , in = 60 ] ( 2 , -2.5 + 0.0329 )
    ( 0 , -2.5 + 0.0329 ) to[ out = 120 , in = 45 ] ( -1 , -2.5 )
    ;

\draw[ thick ]
    ( 0 , 1 + 2.5 ) arc[ radius = 1 , start angle = 90 , end angle = { 270 + 10 } ]
        to[ out = 10 , in = { 180 - 10 } ] ( 2 - 0.1736 , -0.9848 + 2.5 )
        arc[ radius = 1 , start angle = { -90 - 10 } , end angle = { -90 + 10 } ]
        to[ out = 10 , in = { 180 - 10 } ] ( 4 - 0.1736 , -0.9848 + 2.5 )
        arc[ radius = 1 , start angle = { -90 - 10 } , end angle = { -90 + 10 } ]
        to[ out = 10 , in = { 180 - 10 } ] ( 6 - 0.1736 , -0.9848 + 2.5 )
        arc[ radius = 1 , start angle = { -90 - 10 } , end angle = { 90 + 10 } ]
        to[ out = { 180 + 10 } , in = -10 ] ( 4 + 0.1736 , 0.9848 + 2.5 )
        arc[ radius = 1 , start angle = { 90 - 10 } , end angle = { 90 + 10 } ]
        to[ out = { 180 + 10 } , in = -10 ] ( 2 + 0.1736 , 0.9848 + 2.5 )
        arc[ radius = 1 , start angle = { 90 - 10 } , end angle = { 90 + 10 } ]
        to[ out = { 180 + 10 } , in = -10 ] ( 0 + 0.1736 , 0.9848 + 2.5 )
        arc[ radius = 1 , start angle = { 90 - 10 } , end angle = 90 ]
    ( 0 , 1 ) arc[ radius = 1 , start angle = 90 , end angle = { 270 + 10 } ]
        to[ out = 10 , in = { 180 - 10 } ] ( 2 - 0.1736 , -0.9848 )
        arc[ radius = 1 , start angle = { -90 - 10 } , end angle = { -90 + 10 } ]
        to[ out = 10 , in = { 180 - 10 } ] ( 4 - 0.1736 , -0.9848 )
        arc[ radius = 1 , start angle = { -90 - 10 } , end angle = { -90 + 10 } ]
        to[ out = 10 , in = { 180 - 10 } ] ( 6 - 0.1736 , -0.9848 )
        arc[ radius = 1 , start angle = { -90 - 10 } , end angle = { 90 + 10 } ]
        to[ out = { 180 + 10 } , in = -10 ] ( 4 + 0.1736 , 0.9848 )
        arc[ radius = 1 , start angle = { 90 - 10 } , end angle = { 90 + 10 } ]
        to[ out = { 180 + 10 } , in = -10 ] ( 2 + 0.1736 , 0.9848 )
        arc[ radius = 1 , start angle = { 90 - 10 } , end angle = { 90 + 10 } ]
        to[ out = { 180 + 10 } , in = -10 ] ( 0 + 0.1736 , 0.9848 )
        arc[ radius = 1 , start angle = { 90 - 10 } , end angle = 90 ]
    ( 0 , 1 - 2.5 ) arc[ radius = 1 , start angle = 90 , end angle = { 270 + 10 } ]
        to[ out = 10 , in = { 180 - 10 } ] ( 2 - 0.1736 , -0.9848 - 2.5 )
        arc[ radius = 1 , start angle = { -90 - 10 } , end angle = { -90 + 10 } ]
        to[ out = 10 , in = { 180 - 10 } ] ( 4 - 0.1736 , -0.9848 - 2.5 )
        arc[ radius = 1 , start angle = { -90 - 10 } , end angle = { -90 + 10 } ]
        to[ out = 10 , in = { 180 - 10 } ] ( 6 - 0.1736 , -0.9848 - 2.5 )
        arc[ radius = 1 , start angle = { -90 - 10 } , end angle = { 90 + 10 } ]
        to[ out = { 180 + 10 } , in = -10 ] ( 4 + 0.1736 , 0.9848 - 2.5 )
        arc[ radius = 1 , start angle = { 90 - 10 } , end angle = { 90 + 10 } ]
        to[ out = { 180 + 10 } , in = -10 ] ( 2 + 0.1736 , 0.9848 - 2.5 )
        arc[ radius = 1 , start angle = { 90 - 10 } , end angle = { 90 + 10 } ]
        to[ out = { 180 + 10 } , in = -10 ] ( 0 + 0.1736 , 0.9848 - 2.5 )
        arc[ radius = 1 , start angle = { 90 - 10 } , end angle = 90 ] ;

\fill[ Goldenrod , opacity = 0.2 ]
    ( 12 , 0 ) circle[ radius = 4 ] ;

\fill[ White ]
    ( 8 , 0 ) to[ out = 45 , in = 135 ] ( 16 , 0 )
        to[ out = 225 , in = 315 ] ( 8 , 0 ) ;
\fill[ SeaGreen , opacity = 0.2 ]
    ( 8 , 0 ) to[ out = 45 , in = 135 ] ( 16 , 0 )
        to[ out = 225 , in = 315 ] ( 8 , 0 ) ;

\fill[ White ]
    ( 12 , 0 ) circle[ x radius = 4 , y radius = 0.75 ] ;
\fill[ Cyan , opacity = 0.2 ]
    ( 12 , 0 ) circle[ x radius = 4 , y radius = 0.75 ] ;

\draw[ thick ]
    ( 12 , 0 ) circle[ radius = 4 ]
    ( 8 , 0 ) arc[ x radius = 4 , y radius = 0.75 , start angle = 180 , end angle = 360 ]
    ( 12 , -4 ) arc[ x radius = 0.75 , y radius = 4 , start angle = -90 , end angle = 90 ] ;

\draw[ thick , dashed , opacity = 0.7 ]
    ( 16 , 0 ) arc[ x radius = 4 , y radius = 0.75 , start angle = 0 , end angle = 180 ]
    ( 12 , 4 ) arc[ x radius = 0.75 , y radius = 4 , start angle = 90 , end angle = 270 ] ;

\draw[ thick , dotted , RoyalPurple ]
    ( 12 - 1.75 , 2.75 ) node {$ \bullet $} node[ anchor = south ] {$ \omega $}
        to[ out = -45 , in = 90 ] node[ anchor = south east ] {$ > D $} ( 10.9 , 0 ) ;
\draw[ RoyalPurple ]
    ( 10.9 , 0.2 ) to ( 10.7 , 0.2 ) to ( 10.7 , 0 ) ;

\draw
    ( 5 , 2.5 ) node {$ X $}
    ( 5 , 0 ) node {$ Y $}
    ( 5 , -2.5 ) node {$ S $}
    ( 14 , 0 ) node {$ \widetilde{ p } \of{ \CC \of{ S } } $}
    ( 13.5 , -2.5 ) node {$ \CC \of{ X } $} ;

\end{tikzpicture}
\caption{The only simple curves on $ X $ whose projections to $ S $ admit simple elevations to $ Y $ are those in the $ D $-neighborhood of $ \widetilde{ p } \of{ \CC \of{ S } } $ in $ \CC \of{ X } $, which is drawn above in green. For example, the simple curve $ \omega $ drawn on $ X $ has distance more than $ D $ to $ \widetilde{ p } \of{ \CC \of{ S } } $, and so $ p \of{ \omega } $ has no simple elevations to $ Y $.} \label{fig:thm-9}
\end{figure}

We state an example of this phenomenon below, using slightly stronger hypotheses on the covers $ p $ and $ q $ to simplify the statement. Note that the hypothesis below on the subgroups of $ \pi_{ 1 } \of{ S } $ corresponding to the covers $ p $ and $ q $ is satisfied precisely when there is no intermediate cover $ W \to S $ that is covered by both $ X $ and $ Y $. When there are no intermediate covers, $ S $ itself plays the role of each $ W_{ g } $. The constant $ D $ stated therein is explicitly computable using the details in \autoref{sec:effectivizingTang} and the main theorems of \cite{APT22}. See \autoref{fig:thm-9} for a visual description of this result.

\begin{theorem*}\label{thm:effectiveapp}

Suppose that (conjugates of) the subgroups of $ \pi_{ 1 } \of{ S } $ corresponding to the covers $ p $ and $ q $ jointly generate $ \pi_{ 1 } \of{ S } $. Then there is an explicitly computable $ D $, depending only on $ \abs{ \chi \of{ S } } $, $ \deg \of{ p } $, and $ \deg \of{ q } $, so that if $ \gamma $ is any simple closed curve on $ X $ satisfying $ d_{ X } \of{ \gamma , \widetilde{ p } \of{ \CC \of{ S } } } \geq D $, then $ p \of{ \gamma } $ does not admit a simple elevation to $ Y $.

\end{theorem*}

\section{Unmarked simple trace spectral rigidity}
\label{sec:SLnTrace}

The goal of this section is to apply the proof of \autoref{thm:effective} to prove \autoref{thm:SLnTrace}, which we restate below for the reader's convenience.

\SLnTrace*

We will prove an even stronger statement: if $ p $ and $ q $ are not isomorphic, then for $ \K = \R $ and $ \C $, the set of linear representations $ \pi_{ 1 } \of{ S } \to \SL_{ N } \of{ \K } $ whose restrictions to $ \pi_{ 1 } \of{ X } $ and $ \pi_{ 1 } \of{ Y } $ are simple trace isospectral has null Lebesgue measure in $ \hom \of{ \pi_{ 1 } \of{ S } , \SL_{ N } \of{ \K } } $, and its complement is dense. Thus, the term ``generic'' in the above statement of \autoref{thm:SLnTrace} can be interpreted both measure-theoretically and topologically.

In what follows, we will denote by $ R_{ \K } \of{ S , N } $ the \emph{representation variety} $ \hom \of{ \pi_{ 1 } \of{ S } , \SL_{ N } \of{ \K } } $. We will require the following theorem of Rapinchuk--Benyash-Krivetz--Chernousov.

\begin{theorem*}[{\cite[Theorem~3]{RBC96}}] \label{thm:irreducibleRepresentationVariety}

The representation variety $ R_{ \K } \of{ S , N } $ is an irreducible variety.

\end{theorem*}

We will also require the following basic arithmetic lemma.

\begin{lemma*} \label{lem:annoyingAlgebra}

Given complex numbers $ x , y , z \in \C $, if $ x \neq y $, $ z \neq 0 $, and $ x + \frac{ 1 }{ z } = y + z $, then $ 2 x + \frac{ 1 }{ z^{ 2 } } \neq 2 y + z^{ 2 } $.

\begin{proof}

If $ z = \pm 1 $, then the equation $ x + \frac{ 1 }{ z } = y + z $ implies that $ x = y $, and so this contradiction implies that $ z \neq \pm 1 $. Note that 
\begin{align*}
2 x + \frac{ 1 }{ z^{ 2 } } & = 2 \of{ x + \frac{ 1 }{ z } } - \frac{ 2 }{ z } + \frac{ 1 }{ z^{ 2 } } = 2 \of{ y + z } - \frac{ 2 }{ z } + \frac{ 1 }{ z^{ 2 } } = 2 y + z^{ 2 } - z^{ 2 } + 2 z - \frac{ 2 }{ z } + \frac{ 1 }{ z^{ 2 } } \\
& = 2 y + z^{ 2 } - \frac{ z^{ 4 } - 2 z^{ 3 } + 2 z - 1 }{ z^{ 2 } } = 2 y + z^{ 2 } - \frac{ \of{ z - 1 }^{ 3 } \of{ z + 1 } }{ z^{ 2 } } .
\end{align*}
Since $ z \neq \pm 1 $, $ \frac{ \of{ z - 1 }^{ 3 } \of{ z + 1 } }{ z^{ 2 } } $ is non-zero, and so $ 2 x + \frac{ 1 }{ z^{ 2 } } \neq 2 y + z^{ 2 } $, as desired. \qedhere

\end{proof}

\end{lemma*}

\begin{proof}[Proof of \autoref{thm:SLnTrace}]

Fix basepoints $ s \in S $, $ x \in X $, and $ y \in Y $ so that $ s = p \of{ x } = q \of{ y } $, and suppose that $ p \colon X \to S $ and $ q \colon Y \to S $ are non-isomorphic. The above proof of \autoref{thm:effective} implies that there is a curve $ \class{ \gamma } $ whose word length in some chosen finite generating set for $ \pi_{ 1 } \of{ S , s } $ is explicitly bounded in terms of $ \deg \of{ p } $, $ \deg \of{ q } $, and $ \abs{ \chi \of{ S } } $, and with a simple elevation along one of these covers but no simple elevations along the other.

Without loss of generality, assume that $ \class{ \gamma } $ has a simple elevation $ \class{ \alpha } $ along $ p $ to $ X $ but no simple elevations along $ q $ to $ Y $; that is, there is some simple element $ \alpha \in \pi_{ 1 } \of{ X , x } $ so that $ p_{ * } \of{ \alpha } $ is conjugate to $ \gamma^{ m } $ in $ \pi_{ 1 } \of{ S , s } $ for some integer $ 1 \leq m \leq d $. On the other hand, for any simple element $ \beta \in \pi_{ 1 } \of{ Y , y } $, $ q_{ * } \of{ \beta } $ is not conjugate to $ \gamma^{ m } $ or $ \gamma^{ -m } $.

Let $ \Omega \coloneqq \set{ \omega_{ 1 } , \dotsc , \omega_{ r } } \subseteq \pi_{ 1 } \of{ S , s } $ be a set of trace twins of $ \gamma $ for which the commensurability constants are $ m $ and $ 1 $; that is,
\[
\tr \of{ \rho \of{ \gamma^{ m } } }^{ 2 } = \tr \of{ \rho \of{ \omega_{ i } } }^{ 2 },
\]
for all linear representations $ \rho \in R_{ \C } \of{ S , 2 } $ and indices $ 1 \leq i \leq r $. In particular, \cite[Theorem~1.4]{Lei03} implies that $ m \len_{ \mu } \of{ \gamma } = \len_{ \mu } \of{ \omega_{ i } } $ for all $ \mu \in \Teich \of{ S } $ and indices $ 1 \leq i \leq r $.

By \cite[Corollary~3.4]{Lei03}, we may assume that $ \gamma^{ m } $ and $ \omega_{ i } $ are homologous (possibly replacing $ \omega_{ i } $ with $ \omega_{ i }^{ -1 } $ in order to do so), and we may furthermore choose $ \Omega $ to be a maximal such set subject to the additional condition that no two distinct elements are conjugate in $ \pi_{ 1 } \of{ S , s } $. Such a maximal set is finite, since any closed hyperbolic surface has finitely many distinct curves of length at most any finite bound.

\begin{claim} \label{clm:det}

$ \det \of{ \rho \of{ \gamma^{ m } } } = \det \of{ \rho \of{ \omega_{ i } } } $ for any linear representation $ \rho \colon \pi_{ 1 } \of{ S , s } \to \GL_{ n } \of{ \C } $ and any index $ 1 \leq i \leq r $.

\begin{proof}[Proof of \autoref{clm:det}]

The map $ \det \circ \rho \colon \pi_{ 1 } \of{ S , s } \to \C^{ \times } $ is a homomorphism to an abelian group $ \C^{ \times } $, and so it factors through the abelianization $ \Homol_{ 1 } \of{ S ; \Z } $. The desired claim now follows immediately from our choice of orientation of $ \omega_{ i } $ so that $ \gamma^{ m } $ and $ \omega_{ i } $ are homologous. \qedhere

\end{proof}

\end{claim}

Since surface groups and free groups are conjugacy-separable (see \cite{Mar07}, \cite{Dye79}), there is for each index $ 1 \leq i \leq r $ a homomorphism $ \phi_{ i } \colon \pi_{ 1 } \of{ S } \to G_{ i } $ onto a finite group $ G_{ i } $ so that $ \phi_{ i } \of{ \gamma } $ and $ \phi_{ i } \of{ \omega_{ i } } $ are not conjugate in $ G_{ i } $. Note that $ \card{ G_{ i } } $ may, in principle, be explicitly bounded in terms of the word lengths of $ \gamma $ and $ \omega_{ i } $, and that the Milnor-\v{S}varc lemma implies that the word length of $ \omega_{ i } $ is explicitly bounded in terms of the word length of $ \gamma $.

Moreover, basic representation theory implies that there is for each index $ 1 \leq i \leq r $ an irreducible representation $ \psi_{ i } \colon G_{ i } \to \GL_{ N_{ i } } \of{ \C } $ such that
\[
\tr \of{ \psi_{ i } \of{ \phi_{ i } \of{ \gamma^{ m } } } } \neq \tr \of{ \psi_{ i } \of{ \phi_{ i } \of{ \omega_{ i } } } } ,
\]
and whose dimension $ N_{ i } $ is explicitly bounded in terms of the character table for $ G_{ i } $.

Fix $ N \coloneqq \max \set{ 2 , 2 N_{ 1 } + 1 , \dotsc , 2 N_{ r } + 1 } $.

\begin{claim} \label{clm:tr}

For any element $ \delta \in \pi_{ 1 } \of{ S , s } $, either $ \delta $ is conjugate to $ \gamma^{ m } $, or there is some linear representation $ \rho \in R_{ \C } \of{ S , N } $ such that $ \tr \of{ \rho \of{ \gamma^{ m } } } \neq \tr \of{ \rho \of{ \delta } } $.

\begin{proof}[Proof of \autoref{clm:tr}]

Suppose that $ \delta $ is not conjugate to $ \gamma^{ m } $. There are three cases; either $ \gamma $ and $ \delta $ are trace twins with commensurability constants $ m $ and $ 1 $ and $ \gamma^{ m } $ and $ \delta $ are homologous, or they are trace twins with these commensurability constants and $ \gamma^{ m } $ and $ \delta^{ -1 } $ are homologous, or they are not trace twins with these commensurability constants. Suppose first that $ \gamma $ and $ \delta $ are trace twins with commensurability constants $ m $ and $ 1 $; that is,
\[
\tr \of{ \rho \of{ \gamma^{ m } } }^{ 2 } = \tr \of{ \rho \of{ \delta } }^{ 2 }
\]
for all linear representations $ \rho \in R_{ \C } \of{ S , 2 } $, and that $ \gamma^{ m } $ and $ \delta $ are homologous. The maximality of $ \Omega $ implies that $ \delta $ is conjugate to $ \omega_{ i } $ for some index $ 1 \leq i \leq r $. Consider the representation $ \rho \in R_{ \C } \of{ S , N } $ defined by the formula
\begin{equation} \label{eqn:representationFormula}
\rho \of{ \eta } \coloneqq \psi_{ i } \of{ \phi_{ i } \of{ \eta } } \oplus \begin{pmatrix} \frac{ 1 }{ \det \of{ \psi_{ i } \of{ \phi_{ i } \of{ \eta } } } } \end{pmatrix} \oplus \vect{ I }_{ N - N_{ i } - 1 } ;
\end{equation}
that is, $ \rho \of{ \eta } $ is block diagonal, with an $ N_{ i } \times N_{ i } $ block given by $ \psi_{ i } \of{ \phi_{ i } \of{ \eta } } $, a $ 1 \times 1 $ block with entry $ \frac{ 1 }{ \det \of{ \psi_{ i } \of{ \phi_{ i } \of{ \eta } } } } $, and an $ \of{ N - N_{ i } - 1 } \times \of{ N - N_{ i } - 1 } $ identity block. Note that
\begin{align*}
\tr \of{ \rho \of{ \gamma^{ m } } } & = \tr \of{ \psi_{ i } \of{ \phi_{ i } \of{ \gamma^{ m } } } } + \frac{ 1 }{ \det \of{ \psi_{ i } \of{ \phi_{ i } \of{ \eta } } } } + N - N_{ i } - 1 \\
& \neq \tr \of{ \psi_{ i } \of{ \phi_{ i } \of{ \delta } } } + \frac{ 1 }{ \det \of{ \psi_{ i } \of{ \phi_{ i } \of{ \delta } } } } + N - N_{ i } - 1 = \tr \of{ \rho \of{ \delta } } .
\end{align*}

Now suppose that $ \gamma $ and $ \delta $ are trace twins with commensurability constants $ m $ and $ 1 $, but that $ \gamma^{ m } $ and $ \delta^{ -1 } $ are homologous. The maximality of $ \Omega $ implies that $ \delta $ is conjugate to $ \omega_{ i }^{ -1 } $ for some index $ 1 \leq i \leq r $. Consider the representation $ \rho , \rho' \in R_{ \C } \of{ S , N } $, where $ \rho $ is defined by \autoref{eqn:representationFormula}, and $ \rho' $ is defined by the formula
\[
\rho' \of{ \eta } \coloneqq \psi_{ i } \of{ \phi_{ i } \of{ \eta } } \oplus \psi_{ i } \of{ \phi_{ i } \of{ \eta } } \oplus \begin{pmatrix} \frac{ 1 }{ \det \of{ \psi_{ i } \of{ \phi_{ i } \of{ \eta } } }^{ 2 } } \end{pmatrix} \oplus \vect{ I }_{ N - 2 N_{ i } - 1 } .
\]
Therefore either $ \tr \of{ \rho \of{ \gamma^{ m } } } \neq \tr \of{ \rho \of{ \delta } } $ or $ \tr \of{ \rho \of{ \gamma^{ m } } } = \tr \of{ \rho \of{ \delta } } $, in which case \autoref{lem:annoyingAlgebra} implies that $ \tr \of{ \rho' \of{ \gamma^{ m } } } \neq \tr \of{ \rho' \of{ \delta } } $.
Either way we have a representation with the desired property.

Finally, suppose that $ \gamma $ and $ \delta $ are not trace twins with commensurability constants $ m $ and $ 1 $. Then \cite[Theorem~1.4]{Lei03} implies that $ \class{ \gamma } $ and $ \class{ \delta } $ are not length twins with commensurability constants $ m $ and $ 1 $, and so
\[
\tr \of{ \rho' \of{ \gamma^{ m } } }^{ 2 } \neq \tr \of{ \rho' \of{ \delta } }^{ 2 }
\]
for some linear representation $ \rho' \colon \pi_{ 1 } \of{ S , s } \into \SL_{ 2 } \of{ \R } $. Such a linear representation $ \rho' $ can be taken to be a lift to $ \SL_{ 2 } \of{ \R } $ of any monodromy representation of a metric in which $ \gamma^{ m } $ and $ \delta $ have different lengths. In particular, $ \tr \of{ \rho' \of{ \gamma^{ m } } } \neq \tr \of{ \rho' \of{ \delta } } $.

Let $ \iota \colon \SL_{ 2 } \of{ \R } \into \SL_{ N } \of{ \C } $ denote the reducible embedding defined by the formula $ \iota \of{ \vect{ A } } = \vect{ A } \oplus \vect{ I }_{ N - 2 } $, where $ \vect{ I }_{ N - 2 } $ is the $ \of{ n - 2 } \times \of{ n - 2 } $ identity matrix. Note that
\begin{align*}
\tr \of{ \rho \of{ \gamma^{ m } } } & = \tr \of{ \iota \of{ \rho' \of{ \gamma^{ m } } } } = \tr \of{ \rho' \of{ \gamma^{ m } } \oplus \vect{ I }_{ N - 2 } } = \tr \of{ \rho' \of{ \gamma^{ m } } } + n - 2 \\
& \neq \tr \of{ \rho' \of{ \delta } } + n - 2 = \tr \of{ \rho' \of{ \delta } \oplus \vect{ I }_{ n - 2 } } = \tr \of{ \iota \of{ \rho' \of{ \delta } } } = \tr \of{ \rho \of{ \delta } } ,
\end{align*}
where $ \rho \coloneqq \iota \circ \rho' \colon \pi_{ 1 } \of{ S , s } \into \SL_{ N } \of{ \C } $. \qedhere

\end{proof}

\end{claim}

We are now ready to show that $ X $ and $ Y $ are generically simple trace non-isospectral over $ S $ for linear representations $ \pi_{ 1 } \of{ S , s } \to \SL_{ N } \of{ \K } $, where $ \K = \R $ or $ \C $.

Let $ \mathscr{ U }_{ \K } \subseteq R_{ \K } \of{ S , N } $ be the set of linear representations $ \pi_{ 1 } \of{ S , s } \to \SL_{ N } \of{ \C } $ whose restrictions to $ \pi_{ 1 } \of{ X , x } $ and $ \pi_{ 1 } \of{ Y , y } $ along $ p_{ * } $ and $ q_{ * } $, respectively, are simple trace isospectral. We will show that $ \mathscr{ U }_{ \K } $ has null Lebesgue measure in $ R_{ \K } \of{ S , N } $. To that end, for each element $ \delta \in \pi_{ 1 } \of{ S , s } $ not conjugate to $ \gamma^{ m } $, let $ Z_{ \K } \of{ \delta } $ be the zero set of the function $ \tr \of{ \gamma^{ m } } - \tr \of{ \delta } $ on $ R_{ \K } \of{ S , N } $; that is,
\[
Z_{ \K } \of{ \delta } \coloneqq \set{ \rho \in R_{ \K } \of{ S , N } : \tr \of{ \rho \of{ \gamma^{ m } } } = \tr \of{ \rho \of{ \delta } } } .
\]
\autoref{clm:tr} above shows that the regular function $ \tr \of{ \gamma^{ m } } - \tr \of{ \delta } $ on the irreducible variety $ R_{ \C } \of{ S , N } $ is not identically zero. The identity theorem for holomorphic functions implies that it is also not identically zero on the real points $ R_{ \R } \of{ S , N } $.

\begin{claim} \label{clm:positiveCodimension}

For any element $ \delta \in \pi_{ 1 } \of{ S , s } $ not conjugate to $ \gamma^{ m } $, $ Z_{ \K } \of{ \delta } $ is of positive codimension in $ R_{ \K } \of{ S , N } $.

\begin{proof}[Proof of \autoref{clm:positiveCodimension}]

The codimension of $ Z_{ \K } \of{ \delta } $ in $ R_{ \K } \of{ S , N } $ is the height of the principal ideal generated by $ \tr \of{ \gamma^{ m } } - \tr \of{ \delta } $. Since this function is not nilpotent and $ R_{ \K } \of{ S , N } $ is irreducible by \autoref{thm:irreducibleRepresentationVariety}, this height is positive. \qedhere

\end{proof}

\end{claim}

An immediate consequence of \autoref{clm:positiveCodimension} is that for each element $ \delta \in \pi_{ 1 } \of{ S , s } $ not conjugate to $ \gamma^{ m } $, $ Z_{ \K } \of{ \delta } $ is closed and of null Lebesgue measure in $ R_{ \K } \of{ S , N } $, and its complement $ R_{ \K } \of{ S , N } \setminus Z_{ \K } \of{ \delta } $ is open and dense. In particular, as a countable union of null measure sets,
\[
Z_{ \K } \coloneqq \bigcup_{ \delta }{ Z_{ \K } \of{ \delta } },
\]
has null Lebesgue measure, where in the above union, $ \delta $ varies over all elements of $ \pi_{ 1 } \of{ S , s } $ not conjugate to $ \gamma^{ m } $. Moreover, its complement
\[
R_{ \K } \of{ S , N } \setminus Z_{ \K } = R_{ \K } \of{ S , N } \setminus \bigcup_{ \delta }{ Z_{ \K } \of{ \delta } } = \bigcap_{ \delta }{ R_{ \K } \of{ S , N } \setminus Z_{ \K } \of{ \delta } },
\]
is dense in $ R_{ \K } \of{ S , N } $, since $ R_{ \K } \of{ S , N } $ is a Baire space. It therefore suffices to show that $ \mathscr{ U }_\K \subseteq Z_\K $.

To that end, let $ \rho \in \mathscr{ U }_{ \K } $ and $ T \coloneqq \tr \of{ \rho \of{ \gamma^{ m } } } \in \K $. Since
\[
\tr \of{ \of{ \rho \circ p_{ * } } \of{ \alpha } } = \tr \of{ \rho \of{ p_{ * } \of{ \alpha } } } = \tr \of{ \rho \of{ \gamma^{ m } } } = T,
\]
and $ \rho \circ p_{ * } $ and $ \rho \circ q_{ * } $ are simple trace isospectral, there is a simple element $ \beta \in \pi_{ 1 } \of{ S , s } $ so that $ \tr \of{ \of{ \rho \circ q_{ * } } \of{ \beta } } = T $. Let $ \delta \coloneqq q_{ * } \of{ \beta } \in \pi_{ 1 } \of{ S , s } $, and note that
\[
\tr \of{ \rho \of{ \gamma^{ m } } } = \tr \of{ \rho \of{ p_{ * } \of{ \alpha } } } = T = \tr \of{ \rho \of{ q_{ * } \of{ \beta } } } = \tr \of{ \rho \of{ \delta } } .
\]
On the other hand, $ \delta $ cannot be conjugate to $ \gamma^{ m } $, since there are no simple elevations of $ \gamma $ along $ q $ to $ Y $. Thus $ \rho \in Z_{ \K } \of{ \delta } \subseteq Z_{ \K } $. In particular, $ \mathscr{ U }_{ \K } \subseteq Z_{ \K } $, and so $ \mathscr{ U }_{ \K } $ has null Lebesgue measure in $ R_{ \K } \of{ S , N } $, and its complement $ R_{ \K } \of{ S , N } \setminus \mathscr{ U }_{ \K } $ is dense. \qedhere

\setcounter{claim}{0}

\end{proof}

\appendix

\section{Effectivizing an argument of Rivin} \label{sec:Rivin}

In this section, we effectivize the result of Rivin \cite{Riv01} to show that the constants $c_1(\Sigma)$, $c_2(\Sigma)$ are effectively computable. 
We first restate the theorem for the reader's convenience.

\Rivin*

In what follows, we use the strategy of \cite{Riv01} to build the relevant constants recursively from the base case of the sphere with four boundary components $\Sigma_{0,4}$. We outline an effective computation for the constants in the case of the four-holed sphere, and then we note that Rivin's proof for general surfaces allows for a recursive computation of constants from surface to surface. It follows that given explicit constants in the four-holed sphere case as input, $c_{1}(\Sigma)$ and $c_{2}(\Sigma)$ are recursively computable, as desired. 

For the upper bound $c_{2}(\Sigma)$ where $\Sigma$ is some finite type hyperbolic surface, fix a Bers pants decomposition $\mathcal{P}$ on $\Sigma$, and let $\epsilon$ denote the injectivity radius of $\Sigma$. Recall that $\mathcal{P}$ has total geodesic length at most some $B$, where $B$ depends only on the topology of $\Sigma$. Then by basic hyperbolic geometry, there is some constant $\tau = \tau(\epsilon, B)$ so that the minimum length of a geodesic arc connecting two pants curves (or, one pants curve to itself) is at least $\tau$. Then given a simple closed geodesic $\alpha$ with length at most $L$, it follows that 
\[ i(\alpha, \mathcal{P}) \leq L/\tau, \]
and moreover, $\alpha$ can twist about a given pants curve at most $L/B$ times. It therefore follows that every Dehn-Thurston coordinate of $\alpha$ has absolute value bounded above by $L/B + L/\tau$. Thus, $\mathcal{N}(L, \Sigma)$ is at most the number of integer points in a cube in $\mathbb{R}^{\dim(\Teich(\Sigma))}$, which is 
\[ \left( L/\tau + L/B \right)^{6g-6+b+2c}. \] 
One thus obtains a bound on $c_{2}$ of the form 
\[ c_{2}(\Sigma) \leq \left(1/\tau + 1/B \right)^{6g-6+b+2c}. \]

We next turn to the lower bound $c_{1}$. The reader familiar with \cite{Riv01} may notice some subtle differences between Rivin's work and the numerology that follows. 
Since we could not fully reconcile the proof in \cite{Riv01}, we instead opted to describe our understanding of how to implement the basic strategy. 

\begin{figure}[ht]
\centering
\begin{tikzpicture}[ scale = 1.5 ]

\fill[ Orange , opacity = 0.2 ]
    ( 0 , 1.3 ) arc[ x radius = 0.2 , y radius = 0.45 , start angle = 90 , end angle = 270 ]
        arc[ x radius = 0.3 , y radius = 0.4 , start angle = 90 , end angle = -90 ]
        arc[ x radius = 0.2 , y radius = 0.45 , start angle = 90 , end angle = 270 ]
        to[ out = 0 , in = 180 ] ( 1.5 , -0.45 )
        arc[ x radius = 0.2 , y radius = 0.45 , start angle = 270 , end angle = 90 ]
            node[ anchor = south ] {$ E $}
        to[ out = 180 , in = 0 ] ( 0 , 1.3 )
    ( 2.5 , 0.45 ) arc[ x radius = 0.2 , y radius = 0.45 , start angle = 90 , end angle = 270 ]
        to[ out = 0 , in = 180 ] ( 4 , -1.3 )
        arc[ x radius = 0.2 , y radius = 0.45 , start angle = 270 , end angle = 90 ]
        arc[ x radius = 0.3 , y radius = 0.4 , start angle = 270 , end angle = 90 ]
        arc[ x radius = 0.2 , y radius = 0.45 , start angle = 270 , end angle = 90 ]
        to[ out = 180 , in = 0 ] ( 2.5 , 0.45 )
    ( 5 , 1.3 ) arc[ x radius = 0.2 , y radius = 0.45 , start angle = 90 , end angle = 270 ]
        arc[ x radius = 0.3 , y radius = 0.4 , start angle = 90 , end angle = -90 ]
        arc[ x radius = 0.2 , y radius = 0.45 , start angle = 90 , end angle = 270 ]
        to[ out = 0 , in = 180 ] ( 6.5 , -0.45 )
        to[ out = 0 , in = 180 ] ( 8 , -1.3 )
        arc[ x radius = 0.2 , y radius = 0.45 , start angle = 270 , end angle = 90 ]
        arc[ x radius = 0.3 , y radius = 0.4 , start angle = 270 , end angle = 90 ]
        arc[ x radius = 0.2 , y radius = 0.45 , start angle = 270 , end angle = 90 ]
        to[ out = 180 , in = 0 ] ( 6.5 , 0.45 )
        to[ out = 180 , in = 0 ] ( 5 , 1.3 ) ;

\fill[ Orange , opacity = 0.1 ]
    ( 1.5 , 0 ) circle[ x radius = 0.2 , y radius = 0.45 ]
    ( 4 , 0.85 ) circle[ x radius = 0.2 , y radius = 0.45 ]
    ( 4 , -0.85 ) circle[ x radius = 0.2 , y radius = 0.45 ]
    ( 8 , 0.85 ) circle[ x radius = 0.2 , y radius = 0.45 ]
    ( 8 , -0.85 ) circle[ x radius = 0.2 , y radius = 0.45 ] ;

\draw[ thick , RawSienna ]
    ( 1.5 - 0.173 , 0.225 ) to[ out = 180 , in = 345 ] ( 0.1763 , 0.3236 ) ;
\draw[ thick , dashed, RawSienna , opacity = 0.7 ]
    ( 0.1763 , 0.3236 ) to[ out = 15 , in = 185 ] ( 1.5 + 0.1 , 0.3897 ) ;

\draw[ thick , DarkOrchid ]
    ( 1.5 - 0.2 , 0 ) to[ out = 180 , in = 345 ] ( 0.2853 , 0.1236 ) ;
\draw[ thick , dashed , DarkOrchid , opacity = 0.7 ]
    ( 0.2853 , 0.1236 ) to[ out = 15 , in = 185 ] ( 1.5 + 0.173 , 0.225 ) ;

\draw[ thick , MidnightBlue ]
    ( 1.5 - 0.173 , -0.225 ) to[ out = 180 , in = 345 ] ( 0.2853 , -0.1236 ) ;
\draw[ thick , dashed , MidnightBlue , opacity = 0.7 ]
    ( 0.2853 , -0.1236 ) to[ out = 15 , in = 185 ] ( 1.5 + 0.2 , 0 ) ;

\draw[ thick , Green ]
    ( 1.5 - 0.1 , -0.3897 ) to[ out = 180 , in = 345 ] ( 0.1763 , -0.3236 ) ;
\draw[ thick , dashed , Green , opacity = 0.7 ]
    ( 0.1763 , -0.3236 ) to[ out = 15 , in = 185 ] ( 1.5 + 0.173 , -0.225 ) ;

\draw[ thick , Magenta ]
    ( 2.5 - 0.173 , 0.225 ) to[ out = 5 , in = 195 ] ( 4 - 0.1763 , 0.3236 ) ;
\draw[ thick , dashed , Magenta , opacity = 0.7 ]
    ( 4 - 0.1763 , 0.3236 ) to[ out = 165 , in = 355 ] ( 2.5 + 0.1 , 0.3897 ) ;

\draw[ thick , Bittersweet ]
    ( 2.5 - 0.2 , 0 ) to[ out = 5 , in = 195 ] ( 4 - 0.2853 , 0.1236 ) ;
\draw[ thick , dashed , Bittersweet , opacity = 0.7 ]
    ( 4 - 0.2853 , 0.1236 ) to[ out = 165 , in = 355 ] ( 2.5 + 0.173 , 0.225 ) ;

\draw[ thick , Cyan ]
    ( 2.5 - 0.173 , -0.225 ) to[ out = 5 , in = 195 ] ( 4 - 0.2853 , -0.1236 ) ;
\draw[ thick , dashed , Cyan , opacity = 0.7 ]
    ( 4 - 0.2853 , -0.1236 ) to[ out = 165 , in = 355 ] ( 2.5 + 0.2 , 0 ) ;

\draw[ thick , SeaGreen ]
    ( 2.5 - 0.1 , -0.3897 ) to[ out = 5 , in = 195 ] ( 4 - 0.1763 , -0.3236 ) ;
\draw[ thick , dashed , SeaGreen , opacity = 0.7 ]
    ( 4 - 0.1763 , -0.3236 ) to[ out = 165 , in = 355 ] ( 2.5 + 0.173 , -0.225 ) ;

\draw[ thick , RawSienna ]
    ( 5.83 , 0.79 ) to[ out = 180 , in = 55 ] ( 5 + 0.1763 , 0.3236 ) ;
\draw[ thick , dashed , RawSienna , opacity = 0.7 ]
    ( 6.5 + 0.1902 , 0.1391 ) to[ out = 160 , in = 315 ] ( 5.83 , 0.79 )
    ( 5 + 0.1763 , 0.3236 ) to[ out = 15 , in = 165 ] ( 6.5 + 0.1902 , -0.1391 ) ;

\draw[ thick , DarkOrchid ]
    ( 6.2 , 0.522 ) to[ out = 175 , in = 30 ] ( 5 + 0.2853 , 0.1236 ) ;
\draw[ thick , dashed , DarkOrchid , opacity = 0.7 ]
    ( 6.5 + 0.1176 , 0.3641 ) to[ out = 170 , in = 325 ] ( 6.2 , 0.522 )
    ( 5 + 0.2853 , 0.1236 ) to[ out = 0 , in = 160 ] ( 6.5 + 0.1176 , -0.3641 ) ;

\draw[ thick , MidnightBlue ]
    ( 6.5 - 0.1176 , 0.3641 ) to[ out = 185 , in = 15 ] ( 5 + 0.2853 , -0.1236 )
    ( 6.2 , -0.522 ) to[ out = 60 , in = 200 ] ( 6.5 - 0.1176 , -0.3641 ) ;
\draw[ thick , dashed , MidnightBlue , opacity = 0.7 ]
    ( 5 + 0.2853 , -0.1236 ) to[ out = 355 , in = 165 ] ( 6.2 , -0.522 ) ;

\draw[ thick , Green ]
    ( 6.5 - 0.1902 , 0.1391 ) to[ out = 185 , in = 350 ] ( 5 + 0.1763 , -0.3236 )
    ( 5.83 , -0.79 ) to[ out = 105 , in = 195 ] ( 6.5 - 0.1902 , -0.1391 ) ;
\draw[ thick , dashed , Green , opacity = 0.7 ]
    ( 5 + 0.1763 , -0.3236 ) to[ out = 340 , in = 100 ] ( 5.83 , -0.79 ) ;

\draw[ thick , Magenta ]
    ( 6.5 - 0.1902 , 0.1391 ) to[ out = 5 , in = 170 ] ( 8 - 0.1763 , 0.3236 )
    ( 7 , 0.644 ) to[ out = 240 , in = 5 ] ( 6.5 - 0.1176 , 0.3641 ) ;
\draw[ thick , dashed , Magenta , opacity = 0.7 ]
    ( 8 - 0.1763 , 0.3236 ) to[ out = 135 , in = 0 ] ( 7 , 0.644 ) ;

\draw[ thick , Bittersweet ]
    ( 6.5 - 0.1902 , -0.1391 ) to[ out = 15 , in = 175 ] ( 8 - 0.2853 , 0.1236 ) ;
\draw[ thick , dashed , Bittersweet , opacity = 0.7 ]
    ( 8 - 0.2853 , 0.1236 ) to[ out = 170 , in = 350 ] ( 6.5 + 0.1176 , 0.3641 ) ;

\draw[ thick , Cyan ]
    ( 6.5 - 0.1176 , -0.3641 ) to[ out = 20 , in = 180 ] ( 8 - 0.2853 , -0.1236 ) ;
\draw[ thick , dashed , Cyan , opacity = 0.7 ]
    ( 8 - 0.2853 , -0.1236 ) to[ out = 170 , in = 340 ] ( 6.5 + 0.1902 , 0.1391 ) ;

\draw[ thick , SeaGreen ]
    ( 7 , -0.644 ) to[ out = 55 , in = 185 ] ( 8 - 0.1763 , -0.3236 ) ;
\draw[ thick , dashed , SeaGreen , opacity = 0.7 ]
    ( 6.5 + 0.1176 , -0.3641 ) to[ out = 340 , in = 115 ] ( 7 , -0.644 )
    ( 8 - 0.1763 , -0.3236 ) to[ out = 175 , in = 350 ] ( 6.5 + 0.1902 , -0.1391 ) ;

\draw[ thick ]
    ( 0 , 1.3 ) arc[ x radius = 0.2 , y radius = 0.45 , start angle = 90 , end angle = 270 ]
        arc[ x radius = 0.3 , y radius = 0.4 , start angle = 90 , end angle = -90 ]
        arc[ x radius = 0.2 , y radius = 0.45 , start angle = 90 , end angle = 270 ]
        to[ out = 0 , in = 180 ] ( 1.5 , -0.45 )
        arc[ x radius = 0.2 , y radius = 0.45 , start angle = 270 , end angle = 90 ]
            node[ anchor = south ] {$ E $}
        to[ out = 180 , in = 0 ] ( 0 , 1.3 )
    ( 2.5 , 0.45 ) node[ anchor = south ] {$ E $}
        arc[ x radius = 0.2 , y radius = 0.45 , start angle = 90 , end angle = 270 ]
        to[ out = 0 , in = 180 ] ( 4 , -1.3 )
        arc[ x radius = 0.2 , y radius = 0.45 , start angle = 270 , end angle = 90 ]
        arc[ x radius = 0.3 , y radius = 0.4 , start angle = 270 , end angle = 90 ]
        arc[ x radius = 0.2 , y radius = 0.45 , start angle = 270 , end angle = 90 ]
        to[ out = 180 , in = 0 ] ( 2.5 , 0.45 )
    ( 5 , 1.3 ) arc[ x radius = 0.2 , y radius = 0.45 , start angle = 90 , end angle = 270 ]
        arc[ x radius = 0.3 , y radius = 0.4 , start angle = 90 , end angle = -90 ]
        arc[ x radius = 0.2 , y radius = 0.45 , start angle = 90 , end angle = 270 ]
        to[ out = 0 , in = 180 ] ( 6.5 , -0.45 )
        to[ out = 0 , in = 180 ] ( 8 , -1.3 )
        arc[ x radius = 0.2 , y radius = 0.45 , start angle = 270 , end angle = 90 ]
        arc[ x radius = 0.3 , y radius = 0.4 , start angle = 270 , end angle = 90 ]
        arc[ x radius = 0.2 , y radius = 0.45 , start angle = 270 , end angle = 90 ]
        to[ out = 180 , in = 0 ] ( 6.5 , 0.45 )
        to[ out = 180 , in = 0 ] ( 5 , 1.3 )
    ( 6.5 , 0.45 ) node[ anchor = south ] {$ E $}
        arc[ x radius = 0.2 , y radius = 0.45 , start angle = 90 , end angle = 270 ] ;

\draw[ thick , dashed , opacity = 0.7 ]
    ( 0 , 1.3 ) arc[ x radius = 0.2 , y radius = 0.45 , start angle = 90 , end angle = -90 ]
    ( 0 , -1.3 ) arc[ x radius = 0.2 , y radius = 0.45 , start angle = -90 , end angle = 90 ]
    ( 1.5 , -0.45 ) arc[ x radius = 0.2 , y radius = 0.45 , start angle = -90 , end angle = 90 ]
    ( 2.5 , 0.45 ) arc[ x radius = 0.2 , y radius = 0.45 , start angle = 90 , end angle = -90 ]
    ( 4 , 1.3 ) arc[ x radius = 0.2 , y radius = 0.45 , start angle = 90 , end angle = -90 ]
    ( 4 , -1.3 ) arc[ x radius = 0.2 , y radius = 0.45 , start angle = -90 , end angle = 90 ]
    ( 5 , 1.3 ) arc[ x radius = 0.2 , y radius = 0.45 , start angle = 90 , end angle = -90 ]
    ( 5 , -1.3 ) arc[ x radius = 0.2 , y radius = 0.45 , start angle = -90 , end angle = 90 ]
    ( 6.5 , -0.45 ) arc[ x radius = 0.2 , y radius = 0.45 , start angle = -90 , end angle = 90 ]
    ( 8 , 1.3 ) arc[ x radius = 0.2 , y radius = 0.45 , start angle = 90 , end angle = -90 ]
    ( 8 , -1.3 ) arc[ x radius = 0.2 , y radius = 0.45 , start angle = -90 , end angle = 90 ] ;

\draw
    ( 2 , 0 ) node {$ + $}
    ( 4.5 , 0 ) node {$ = $}
    ( 0 + 0.2 , -0.85 ) node[ anchor = west ] {$ P_{ 1 } $}
    ( 5 + 0.2 , -0.85 ) node[ anchor = west ] {$ P_{ 1 } $}
    ( 4 - 0.2 , -0.85 ) node[ anchor = east ] {$ P_{ 2 } $}
    ( 8 - 0.2 , -0.85 ) node[ anchor = east ] {$ P_{ 2 } $}
    ( 0 + 0.3 , 0 ) node[ anchor = east ] {$ \gamma_{ 1 } $}
    ( 4 - 0.3 , 0 ) node[ anchor = west ] {$ \gamma_{ 2 } $}
    ( 8 - 0.3 , 0 ) node[ anchor = west ] {$ \gamma_{ p / 2 k } $} ;

\end{tikzpicture}
\caption{The construction of $ \gamma_{ p / 2 k } $ (with $ p = 3 $ and $ k = 4 $)} \label{fig:twistyPants}
\end{figure}
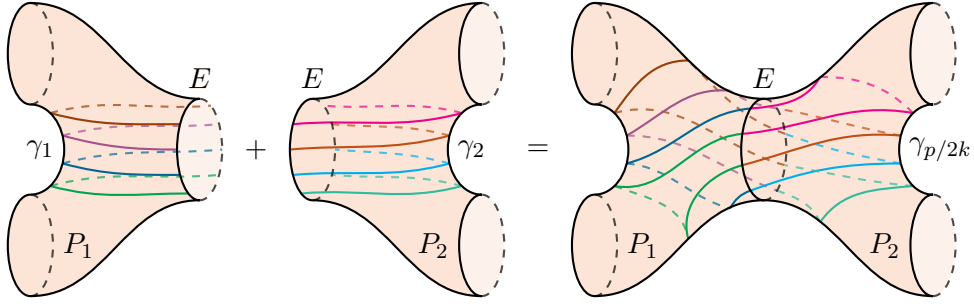

In what follows, we use $E$ to denote a simple separating curve cutting $\Sigma_{0,4}$ into two pair of pants $\mathcal{P}_1$, $\mathcal{P}_2$ glued along $E$.
Note that any simple closed curve on $\Sigma_{0,4}$ is either a power of a boundary component or is obtained by the following procedure (an example of which is drawn in \autoref{fig:twistyPants}): 
\begin{itemize}
\item Fix $k\in\N$ and $p\in\mathbb{Z}$,
\item For $i\in\{1,2\}$, let $\gamma_i$ be $k$ copies of the seam $\eta$ connecting $E$ to $E$ and separating the other boundary components in $\mathcal{P}_i$, where we ask that each of the $k$ seams has a distinct endpoint on $E$ and that the seams are pairwise disjoint,
\item Glue $\gamma_1$ to $\gamma_2$ using a $p/2k$-twist along $E$, where this twist is clockwise if $p/2k>0$ and this twist is counter-clockwise if $p/2k<0$.
\end{itemize}
We will denote a curve obtained by this procedure by $\gamma_{p/2k}$ and note that this procedure gives a simple closed curve precisely when $\gcd(p,2k)=1$.
Recall that for any fixed $k\in\N$, $\phi(2k)$ is the number of $1\le p\le 2k$ such that $\gcd(p,2k)=1$ and hence for any $N\in\N$
$$\left|\left\{\gamma_{p/2k}\mid \left|\frac{p}{2k}\right|<N\right\}\right|=\begin{cases}
4N-1,&k=1,\\
2\phi(2k)N,&k\neq 1.
\end{cases}$$
Indeed the $k\neq 1$ case follows because there are precisely $2\phi(2k)$ possible fractional twists (twists with $-1<p/2k<1$), $\phi(2k)$ positive and $\phi(2k)$ negative, and therefore there are $2\phi(2k)N$ curves $\gamma_{p/2k}$ such that $-N<p/2k<N$.
The $k=1$ case differs because only half integer twists are allowed and so there are $4N-1$ half integers with absolute value less than $N$.
Note that in either case, the size of this set if bounded below by $2\phi(2k)N$.

Now let $m$ by a hyperbolic metric on $\Sigma_{0,4}$ and note that the length of $\gamma_{p/2k}$ is bounded above by
\begin{equation}\label{eqn:gammapqlength}
\ell(\gamma_{p/2k})\le 2k\ell(\eta)+|p|\ell(E),
\end{equation}
where $\eta$ is the geodesic representative of the aforementioned seam.
Let $\mathcal{N}_{0,4}(L,m)$ denote the size of the set of isotopy classes of non-boundary parallel simple closed geodesics of length at most $L$.
Suppose moreover that we require that $m$ is such that $\ell(\eta)\le\ell(E)$. Then it follows from \autoref{eqn:gammapqlength} that
$$\ell(\gamma_{p/2k})\le (2k+|p|)\ell(E).$$
Fix some $L\in\mathbb{R}$, then if $N$ is any natural number for which $N+1\le \frac{L}{2k\ell(E)}$, it follows that if $|p/2k|=|p|/2k< N$ then
$$\ell(\gamma_{p/2k})\le (2k+|p|)\ell(E)< 2k(N+1)\ell(E)\le L.$$
Moreover, by the calculation from the previous paragraph, there are at least
$$\left\lfloor \frac{L}{2k\ell(E)} -1 \right\rfloor \cdot 2\phi(2k),$$
such curves.
Note that the maximum $k$ for which $\ell(\gamma_{p/2k})\le \ell$ for all $-1<p/2k<1$ is $\left\lfloor\frac{L}{4\ell(E)}\right\rfloor$.
Combining the above, it follows that
$$\mathcal{N}_{0,4}(L,m)=\{\gamma\mid \gamma\text{ is a simple closed curve on }\Sigma_{0,4}\}\ge A(L)=\sum_{k=1}^{\left\lfloor\frac{L}{4\ell(E)}\right\rfloor}\left\lfloor \frac{L}{2k\ell(E)} -1 \right\rfloor \cdot 2\phi(2k).$$
We now assume that $L$ is a multiple of $4\ell(E)$ so that  
$$\mathcal{N}_{0,4}(L,m)\ge A(L)=\sum_{k=1}^{\frac{L}{4\ell(E)}} \left(\frac{L}{2k\ell(E)} -1\right) \cdot 2\phi(2k),$$
then importantly $A(L)$ is asymptotic to $\frac{3\pi ^2}{2\ell(E)^2}L^2$.
One can show this rigorously using some basic number theoretic calculations however, because we want an effective theorem, we will not prove this asymptotic directly but rather content ourselves to find a lower bound on $A(L)$ which is quadratic in $L$ for large enough $L$.

For ease of notation define $x=\frac{L}{4\ell(E)}$, then decomposing $A(L)$ as $$A(L)=A^+(L)-A^-(L)=4x\sum_{k=1}^{x}\frac{\phi(2k)}{k}-2\sum_{k=1}^{x}\phi(2k),$$
then it is straightforward that
$$A^-(L)=2\sum_{k=1}^{x}\phi(2k)\le 4\sum_{k=1}^{x}\phi(k)=4\Phi(x).$$
By the M\"obius inversion formula, 
$$\Phi(x)=\sum_{k=1}^{x}\phi(k)=\frac{1}{2}\sum_{k=1}^x\mu(k)\left\lfloor\frac{x}{k}\right\rfloor\left(\left\lfloor\frac{x}{k}\right\rfloor+1\right)=\frac{1}{2}\sum_{k=1}^x\mu(k)\left\lfloor\frac{x}{k}\right\rfloor^2+\frac{1}{2}\sum_{k=1}^x\mu(k)\left\lfloor\frac{x}{k}\right\rfloor.$$
Note that
$$\frac{1}{2}\sum_{k=1}^x\mu(k)\left\lfloor\frac{x}{k}\right\rfloor\le \frac{1}{2}\sum_{k=1}^x\frac{x}{k}\le \frac{1}{2}x(\ln(x)+1).$$
Similarly,
\begin{align*}
\frac{1}{2}\sum_{k=1}^x\mu(k)\left\lfloor\frac{x}{k}\right\rfloor^2&=\frac{1}{2}\sum_{k=1}^x\mu(k)\left(\left\lfloor\frac{x}{k}\right\rfloor^2-\left(\frac{x}{k}\right)^2\right)+\frac{1}{2}\sum_{k=1}^x\mu(k)\left(\frac{x}{k}\right)^2,\\
&=\frac{1}{2}\sum_{k=1}^x\mu(k)\left(\left\lfloor\frac{x}{k}\right\rfloor-\frac{x}{k}\right)\left(\left\lfloor\frac{x}{k}\right\rfloor+\frac{x}{k}\right)+\frac{x^2}{2}\sum_{k=1}^x\frac{\mu(k)}{k^2},\\
&\le\sum_{k=1}^x\frac{x}{k}+\frac{x^2}{2}\left(\frac{6}{\pi^2}-\sum_{k=x+1}^\infty\frac{\mu(k)}{k^2}\right),\\
&\le x(\ln(x)+1)+\frac{3}{\pi^2}x^2+\frac{x^2}{2}\sum_{k=x+1}^\infty\frac{1}{k^2},\\
&\le x(\ln(x)+1)+\frac{3}{\pi^2}x^2+\frac{x^2}{2}\frac{1}{x-1}.
\end{align*}
Therefore, provided $x\ge 2$ to simplify the last term in the previous inequality, we find that
$$A^-(L)\le \frac{12}{\pi^2}x^2+6x\ln(x)+8x+4.$$
For $A^+(L)$, we similarly see that
$$A^+(L)=4x\sum_{k=1}^x\frac{\phi(2k)}{k}\ge4x\sum_{k=1}^x\frac{\phi(k)}{k}=4x\Psi(x).$$
Again by the M\"obius inversion formula
$$\Psi(x)=\sum_{k=1}^x\frac{\phi(k)}{k}=\sum_{k=1}^x\left\lfloor\frac{x}{k}\right\rfloor\frac{\mu(k)}{k}.$$
In particular, similar to the case of $\Phi(x)$ we conclude that
\begin{align*}
\Psi(x)&=\sum_{k=1}^x\left(\left\lfloor\frac{x}{k}\right\rfloor-\frac{x}{k}\right)\frac{\mu(k)}{k}+\sum_{k=1}^x\frac{x}{k}\frac{\mu(k)}{k},\\
&\ge -\sum_{k=1}^x\frac{1}{k}+x\sum_{k=1}^x\frac{\mu(k)}{k^2},\\
&\ge -\ln(x)+x\left(\frac{6}{\pi^2}-\sum_{k=x+1}^\infty\frac{\mu(k)}{k^2}\right),\\
&\ge -\ln(x)+x\left(\frac{6}{\pi^2}-\sum_{k=x+1}^\infty\frac{1}{k^2}\right),\\
&\ge-\ln(x)+x\left(\frac{6}{\pi^2}-\frac{1}{x-1}\right),\\
&\ge -\ln(x)+\frac{6}{\pi^2}x-1.
\end{align*}
Therefore we find that 
$$A^+(L)=4x\Psi(x)\ge \frac{24}{\pi^2}x^2-4x\ln(x)-4x,$$
and consequently
$$A(L)\ge A^+(L)-A^-(L)\ge \frac{12}{\pi^2}x^2-10x\ln(x)-12x-4,$$
which, when remembering that $x=\frac{L}{4\ell(E)}$, has leading term of the form $\frac{3\pi^2}{4\ell(E)}L^2.$
If $x\ge 30$, then $-2x\ln(x)-12x-4\ge -\frac{6}{\pi^2}x^2$ and consequently if we assume that $L\ge L_0=120\ell(E)$, then we find that 
$$A(L)\ge \frac{3}{8\pi^2\ell(E)^2}L^2=c_1(\Sigma_{0,4})L^2,$$
as required.

\section{Effectivizing an argument of Tang} \label{sec:effectivizingTang}

The goal of this appendix is to make clear that the bound in \autoref{thm:Tang} of Tang, which we restate below for the reader's convenience, can be made explicitly computable in the sense of \autoref{sec:effective}.

\Tang*

The \emph{radius} $ \rad \of{ U } $ of a non-empty finite subset $ U \subseteq X $ of a Gromov hyperbolic metric space $ X $ is defined by
\[
\rad \of{ U } \coloneqq \min \set{ r \in \clopint{ 0 }{ \infty } : U \subseteq \mD_{ r } \of{ x } \textrm{ for some } x \in X } .
\]
We say that a point $ x \in X $ is a \emph{circumcenter} of $ U $ if $ U $ is contained in a closed disk of radius $ \rad \of{ U } $ and center $ x $.

In the remainder of this section, we will briefly explain how to make \autoref{thm:Tang} effective.

\subsection{Definitions}

In \cite{Tan17}, Tang defines the \emph{hull} of a tuple of simple curves $ \alpha = \of{ \alpha_{ 1 } , \dotsc , \alpha_{ n } } $ on $ X $ to be the set of all simple curves on $ X $ lying on some geodesic segment connecting two curves in $ \alpha $. Given such a tuple $ \alpha $, fix $ \vect{ t } = \of{ t_{ 1 } , \dotsc , t_{ n } } \in \R_{ \geq 0 }^{ n } $, and let $ \vect{ t } \cdot \alpha $ denote the formal sum
\[
\vect{ t } \cdot \alpha \coloneqq \sum_{ k = 1 }^{ n }{ t_{ k } \alpha_{ k } } .
\]
Tang also defines $ \norm{ \vect{ t } }_{ \alpha } $ by the formula
\[
\norm{ \vect{ t } }_{ \alpha } \coloneqq \sqrt{ \sum_{ j , k = 1 }^{ n }{ t_{ j } t_{ k } i \of{ \alpha_{ j } , \alpha_{ k } } } } ,
\]
and the \emph{$ L $-short set} of $ \vect{ t } \cdot \alpha $ by
\[
\short \of{ \vect{ t } \cdot \alpha , L } \coloneqq \set{ \gamma \in \CC \of{ X } : i \of{ \vect{ t } \cdot \alpha , \gamma } \leq L \norm{ \vect{ t } }_{ \alpha } } .
\]

\subsection{The constants in \cite[Lemma~6.4]{Tan17} and \cite[Lemma~5.1]{Tan17}}

In this subsection, we discuss a few lemmas from \cite{Tan17} and the dependencies of various constants that arise therein. We recall that $ p \colon X \to S $ is a finite-degree cover of a closed surface $ S $ of negative Euler characteristic, and we begin with the following lemma.

\begin{lemma*}[{\cite[Lemma~6.4]{Tan17}}] \label{lem:TangShort}

There is a positive number $ L_{ 0 } > 0 $, depending only on $ \abs{ \chi \of{ X } } $, so that
\[
p \of{ \alpha } \in \short \of{ \sum_{ g \in \Deck \of{ p } }{ g \of{ \alpha } } , L_{ 0 } \card{ \Deck \of{ p } } }
\]
for all simple curves $ \alpha $ on $ X $.

\end{lemma*}

By following the proof of \cite[Lemma~6.4]{Tan17}, one may observe that this constant $ L_{ 0 } $ is the same as that which appears in \cite[Lemma~5.1]{Tan17}.

\begin{lemma*}[{\cite[Lemma~5.1]{Tan17}}]

There are some positive numbers $ L_{ 0 } $ and $ k_{ 0 } $, depending only on $ \abs{ \chi \of{ X } } $, so that for all $ L \geq L_{ 0 } $, the $ L $-short set $ \short \of{ \vect{ t } \cdot \alpha , L } $ is non-empty and has diameter
\[
\diam_{ \CC \of{ X } } \of{ \short \of{ \vect{ t } \cdot \alpha , L } } \leq 4 \log_{ 2 } \of{ L } + k_{ 0 } .
\]

\end{lemma*}

This is proved in \cite[Section~8]{Tan17}. In that section, Tang shows that it suffices to choose $ L_{ 0 } \geq \frac{ \sqrt{ 2 } }{ W_{ 0 } } $, where $ W_{ 0 } > 0 $ depends only on $ \abs{ \chi \of{ S } } $ and satisfies the following condition:

\begin{proposition*}[{\cite[Proposition~4.1(3)]{Tan17}}] \label{prop:TangAnnulus}

$ S \of{ \vect{ t } \cdot \alpha } $ contains an essential annulus of width at least $ W_{ 0 } \norm{ \vect{ t } }_{ \alpha } $.

\end{proposition*}

Above, $ S \of{ \vect{ t } \cdot \alpha } $ is the singular flat surface dual to the union of curves in $ \alpha $, where the rectangle corresponding to an intersection point between $ \alpha_{ j } $ and $ \alpha_{ k } $ has side lengths $ t_{ j } $ and $ t_{ k } $. This proposition is proved using the following lemma of Bowditch.

\begin{lemma*}[{\cite[Lemma~5.1]{Bow06}}] \label{lem:BowditchAnnulus}

Let $ \rho $ be a singular Riemannian metric on an orientable closed surface $ \Sigma $ with area $ 1 $, and let $ A \subsetneq \Sigma $ be a finite set. If there is a homeomorphism $ f \colon \clopint{ 0 }{ \infty } \to \clopint{ 0 }{ \infty } $ so that
\[
\area \of{ D } \leq f \of{ \len \of{ \partial D } }
\]
for any region $ D $ of $ \Sigma $ containing at most one point in $ A $, then $ \Sigma \setminus A $ contains an essential annulus of width at least $ \eta > 0 $, where $ \eta $ depends only on $ \abs{ \chi \of{ \Sigma } } $ and $ \card{ A } $.

\end{lemma*}

In our desired context, the set $ A $ of marked points is empty, and so \autoref{prop:TangAnnulus} follows from \autoref{lem:BowditchAnnulus} once one has shown that $ \area \of{ D } \leq 4 \len \of{ \partial D }^{ 2 } $ for all appropriate regions $ D $ in $ S \of{ \vect{ t } \cdot \alpha } $.

\subsection{A weaker version of \cite[Lemma~5.1]{Tan17}}

The second clause to \cite[Lemma~5.1]{Tan17} states that $\mbox{short}(\textbf{t} \cdot \alpha, L)$ has uniformly bounded diameter in $\mathcal{C}(\Sigma)$ \textit{in terms only of $L$} and $\Sigma$ and not on \textbf{t}. We will need this second clause; however, Tang's proof is complicated by the possibility that $\alpha$ does not fill the surface. Since in our context, $\alpha$ will always be the $G$-orbit of some curve in $\mathcal{C}(\Sigma)$, we can simplify matters by proving a basic lemma which guarantees that $G \cdot \gamma$ fills $\Sigma$ once $\gamma$ is at least distance $2$ from $\pi(\gamma)$:

\begin{lemma*} \label{Fills} Suppose $d_{\Sigma}(\gamma, \pi(\gamma)) \geq 2$, then $G \cdot \gamma$ fills $\Sigma$. 
\end{lemma*}

\begin{proof} If $G \cdot \gamma$ does not fill $\Sigma$, then there is some essential subsurface $Y$ in its complement. Restricting the covering map $P$ to $Y$ yields a cover $P|_{Y}: Y \rightarrow Z$ where $Z$ is some non-simply connected subsurface of $S$. Thus, $Z$ supports a (potentially boundary parallel) simple closed curve $\eta$ which is necessarily disjoint from $P(\gamma)$. It follows that $\widetilde{P}(\eta)$ (which is contained in $\widetilde{P}(\mathcal{C}(S))$) is distance $1$ from some elevation of $\gamma$ to $\Sigma$. 
\end{proof}

\subsection{Assembling the proof modulo $ W_{ 0 } $ and Bowditch's lemma}

In this subsection, we will obtain an effective version of \cite[Proposition~6.2]{Tan17}, save for the constant $W_{0}$. First, we effectivize the proof of \cite[Lemma~5.1]{Tan17} in the event that $\alpha$ fills. There, Tang shows that any two curves in $\mbox{short}(\textbf{t} \cdot \alpha, L)$ (for $L> L_{0}$ with $L_{0}$ as above) intersect at most $L/W_{0}$ times. 

We can then use an effective version of Hempel's bound on curve complex distance in terms of intersection number, due to Bowditch \cite{Bow14} (this particular version can be found in \cite{MT22})
\[ d_{\CC \of{\Sigma}}(\alpha, \beta) < 6 + 2 \cdot \frac{\log(2 \cdot i(\alpha, \beta) + 1)}{\log((2g+p-4)/2)}, \]
where $g,p$ denotes the genus and number of punctures of $\Sigma$, respectively. Tang points out that $L$ need only be bigger than $\sqrt{2}/W_{0}$, and Bowditch obtains an effective bound on $W_{0}$ in \cite{Bow14} which is the following bound
\[ W_{0} \geq \frac{1}{8\cdot (2g+p-1) \cdot (2g+p+6)}. \]

 We note that this also uses Tang's isoperimetric inequality mentioned above; Bowditch states the result in terms of some constant $h$ so that the area of a trivial region is at most $h$ times the length of its boundary squared and moreover shows that a factor of $4\sqrt{h}$ should appear in the bound. Thus, the appearance of $8$ in the denominator above is actually $4 \sqrt{4}$, where the $\sqrt{4}$ corresponds to the $4$ in Tang's isoperimetric bound.
 
Lemma 6.4 of \cite{Tan17} above therefore says that $\pi(\alpha)$ lies in $\mbox{short}({\bf{1}} \cdot G\alpha, L_{0}|G|)$, which by the above, has diameter at most 
\[ 12+ 4 \cdot \frac{\log(2 \cdot (8 \sqrt{2}\cdot (2g+p-1) \cdot (2g+p+6))^{2}\cdot |G| + 1)}{\log((2g+p-4)/2)}, \]
which is bounded uniformly above independent of $g,p$.

We would be done if the circumcentre of $G \cdot \alpha$ were also in this short curve set. It need not be, but \cite[Proposition~5.2]{Tan17} asserts the existence of some $k_{1}$ so that for all $L> L_{0}$, $\mbox{short}(\rho, L)$ (which is defined to be the union, taken over all \textbf{t}, of $\mbox{short}(\textbf{t} \cdot \rho, L)$) is Hausdorff distance at most $k_{1}$ from $\mbox{Hull}(\rho)$, for any finite collection of curves $\rho$. Here $k_{1}$ depends only on $L$, the topology of the surface, and the number of curves in $\rho$. In particular when $\rho = G\cdot \alpha$, $k_{1}$ will depend only on $|G|$, the topology of $\Sigma$ (or alternatively, the topology of $S$ and the degree of $P$), and $L_{0}$. 

In \cite[Proposition~5.2]{Tan17}, Tang shows that $k_{1}$ can be taken to be a constant $r = r(\Sigma, |G|, L_{0})$, where $r$ is defined to be the distance between 
\[ \mbox{short}(\mbox{\textbf{t}}' \cdot \alpha', |G|L_{0}/\sqrt{2}), \]
and $[\alpha_{1}, \alpha_{2}]$, where
\begin{enumerate}
    \item $\alpha'= \left\{\alpha_{1}, \alpha_{2} \right\}$, the points in the orbit of $G \cdot \alpha$ maximizing $t_{i}t_{j} i(\alpha_{i}, \alpha_{j})$ over all pairs $(i,j) \in \left\{1,..., |G| \right\}^{2}$, and 
    \item $[\alpha_{1}, \alpha_{2}]$ is any geodesic in the curve graph between $\alpha_{1}, \alpha_{2}$. 
\end{enumerate}

Since curve graphs are $17$-hyperbolic, any two geodesics are at most $92 \cdot 17= 1564$ apart (see for instance \cite[Theorem~1.1]{GS19}), so up to addition of this universal constant, the choice of geodesic will not matter.

\subsection{Effectivizing the Hausdorff distance statement from Bowditch} 

It remains to determine how close the short curve set $\mbox{short}(\left\{\alpha_{1}, \alpha_{2}\right\}, |G|W_{0}^{-1})$ is to a geodesic segment $[\alpha_{1}, \alpha_{2}]$. For this, we use \cite[Proposition~6.2]{Bow06}, which asserts that there is some $q$ and some $h$ so that for all $Q> q$, the set
\[ G_{Q}(\alpha, \beta) = \left\{ \gamma \in \mathcal{C} : i(\gamma, \alpha) \cdot i(\gamma, \beta) < Q \cdot i(\alpha, \beta) \right\}, \]
is within $h$ of $[\alpha, \beta]$. To understand the relevance of this, first consider some $\gamma \in \mbox{short}(\left\{\alpha, \beta\right\}, Q)$. Thus, 
\[ i(\gamma, t_{1} \alpha) + i(\gamma, t_{2} \beta) < Q \cdot \sqrt{t_{1}t_{2}i(\alpha, \beta)}, \]
for some nonnegative $t_{1}, t_{2}$. Squaring both sides yields 
\[ t_{1}t_{2}i(\gamma, \alpha) \cdot i(\gamma, \beta) + T < Q^{2}t_{1}t_{2} i(\alpha, \beta), \]
where $T$ is non-negative. Therefore, it follows that $\gamma \in G_{Q^{2}}(\alpha, \beta)$. 

Bowditch proves \cite[Proposition~6.2]{Bow06} using the ``lines'' he constructs in the proof of hyperbolicity, which he denotes by $\Lambda_{\alpha \beta}$. In \cite[Section~3]{Bow06}, he outlines a series of criteria involving a family of lines such that if any graph admits such a family, it is hyperbolic, and the lines are uniformly quasigeodesic. The criteria appear directly before \cite[Proposition~3.1]{Bow06}.

The proof of \cite[Proposition~3.1]{Bow06} uses an auxiliary set of paths denoted $\pi_{\alpha \beta}$. In \cite[Lemma~3.3]{Bow06}, Bowditch shows that these paths are $(K, 2K)$-quasigeodesic, where $K$ is as in the statement of the three axioms he uses as a criteria for hyperbolicity. He also shows, directly before this lemma, that $\Lambda_{\alpha \beta}$ is Hausdorff distance at most $K$ from $\pi_{\alpha \beta}$. Putting all of this together and using the Morse Lemma, it follows that $\Lambda_{\alpha, \beta}$ is within a distance of $92K^{2}(2K+2 + 17)$ of a geodesic $[\alpha, \beta]$ (again using \cite{GS19}). 
In \cite[Section~4]{Bow06}, Bowditch shows that the curve graph satisfies the three required axioms, for $K= 18D+2$. 

Furthermore, Bowditch shows in \cite{Bow14} that it suffices to choose $D= 20$ for all surfaces. Putting all of this together, we get that $\Lambda_{\alpha \beta}$ is within a neighborhood of $[\alpha, \beta]$ of radius at most 
\begin{align*} 92K^{2}(2K+19) &= 92\cdot (18D+2)^{2} \cdot (36D+ 23), \\
 &= 8,957,643,664. \end{align*}

In the direction of the proof of \cite[Proposition~6.2]{Bow06} that is relevant to us, Bowditch only needs to establish that $G_{Q}(\alpha, \beta)$ is nearby $\Lambda_{\alpha \beta}$, since as mentioned above, we can estimate the distance between $\Lambda_{\alpha \beta}$ and $[\alpha, \beta]$. 

Bowditch's proof of \cite[Proposition~6.2]{Bow06} shows that one may set $q = (W_{0}^{-1})^{2}$, and then that 
\[ G_{Q}(\alpha, \beta) \subset L_{\alpha \beta}(t, \sqrt{Q}),\]
for all $Q > q$, where $L_{\alpha \beta}(t, \sqrt{Q})$ is as defined in \cite[Section~4]{Bow06}. Thus we are interested in the diameter of $L_{\alpha \beta}(t, (W_{0}^{-1})^{2})$, which, again by Bowditch's proof of uniform hyperbolicity \cite{Bow14}, is uniformly bounded above, for instance by $20$.  

\subsection{Putting it all together.} Now, assume we begin with some $\gamma \in \mbox{short}(\left\{\alpha_{1}, \alpha_{2} \right\}, |G|W_{0}^{-1})$. Then 
\[ \gamma \in G_{|G|^{2}(W_{0}^{-1})^{2}}(\alpha, \beta). \]
It follows from the very end of the previous section, and \cite[Lemma~4.2]{Bow06}, that $\gamma$ is at most
\[ 2(W_{0}^{-1})^{3}|G|+ 2,\]
away from $\Lambda_{\alpha \beta}$. Therefore, $\gamma$ is at most 
\[ 2(W_{0}^{-1})^{3}|G|+ 2 + 8,957,643,664,  \]
away from $[\alpha, \beta]$. 

It finally follows, assuming $G \cdot \alpha$ fills $\Sigma$, that $\pi(\alpha)$ is within 
\[ 2(W_{0}^{-1})^{3}|G|+ 2 + 8,957,643,664 + 12+ 4 \cdot \frac{\log(2 \cdot (8 \sqrt{2}\cdot (2g+p-1) \cdot (2g+p+6))^{2}\cdot |G| + 1)}{\log((2g+p-4)/2)} ,   \]
of the circumcentre of $G \cdot \alpha$.

If $G \cdot \alpha$ does not fill $\Sigma$, \autoref{Fills} implies that the radius of $G \cdot \alpha$ is at most $2$ and that $\pi(\alpha)$ is distance at most $2$ away, and so the circumcentre of $G \cdot \alpha$ is at most distance $4$ from $\pi(\alpha)$.  

Since the log term above is bounded above independently of $g,p$, the entire bound is $O(|G|)$ and on the order of a polynomial of degree $6$ in $g$ and $p$. For a less explicit upper bound that is stronger asymptotically, note that $|G|$ is necessarily bounded above by a linear function of $g$ and $p$. Therefore $|G|^{2} \cdot (W_{0}^{-1})^{2}$ grows at most polynomially in $g,p$ of degree at most $6$. The techniques in \cite{Aou13} show how to obtain a bound that is uniform in $g,p$ on the diameter of $L_{\alpha \beta}(t, R)$ whenever $R$ grows slower than some polynomial in $g,p$. 

\printbibliography

\end{document}